\documentclass[12pt]{amsart}
\usepackage{amssymb,hyperref,times}
\usepackage{latexsym}
\usepackage{amsmath}
\usepackage{amssymb}
\usepackage{amssymb, amscd}
\usepackage[dvips]{graphicx}
\usepackage{ifthen}
\usepackage{amsthm}
\input epsf
\usepackage[T1]{fontenc}

\newtheorem{prop}{Proposition}[subsection]
\newtheorem{lemma}{Lemma}[subsection]

\def\bu{{\bf u}}
\def\bv{{\bf v}}
\def\bw{{\bf w}}
\def\ent{{\rm Ent} T}
\def\Th{\Theta}
\def\cal{\mathcal}
\def\al{\alpha}
\def\G{\mathcal G}
\def\ga{\gamma}
\def\ti{\widetilde}
\def\TG{\widetilde{\G}}
\def\TD{\widetilde{\Delta}}
\def\T{\mathcal T}
\def\N{\Bbb N}
\def\La{\b\Lambda}
\def\hn{{\Bbb H}^n}
\def\h3{{\Bbb H}^3}
\def\tupl#1{\langle #1\rangle}
\let\s\mathsf
\let\Cal\mathcal
\let\g\mathfrak
\let\on\curvearrowright
\def\d{\boldsymbol\delta}
\def\l#1{\vskip3pt\noindent\bf Lemma \the\secNum.#1\sl.}
\def\lnk#1{$^!$}

\font\cFont=cmr9%
\def\note#1{\{{\cFont #1}\}}%

\title{Quasiconvexity in the relatively hyperbolic groups}

\author
{Victor Gerasimov and Leonid Potyagailo}
\address{Victor Guerassimov, Departamento de Matem\'atica,
Universidade Federal de Minas Gerais, Av. Antonio Carlos, 6627/
Caixa Postal 702  30161-970 Belo Horizonte, MG, Brasil}
\email{victor@mat.ufmg.br}
\address{ Leonid Potyagailo, UFR de Math\'ematiques, Universit\'e
de Lille 1, 59655 Villeneuve d'Ascq cedex, France}
\email{potyag@math.univ-lille1.fr}

\subjclass[2000]{Primary 20F65, 20F67; Secondary 57M07, 30F40}
\keywords{quasigeodesic, horosphere, horocycle, quasiconvexity,
$\alpha$-distorted map.}

\date{\today}

\begin{document}

\begin{abstract}
 We study different notions of quasiconvexity  for a subgroup $H$ of
 a relatively hyperbolic group $G.$ The first result establishes
 equivalent conditions for $H$ to be relatively quasiconvex. As a
 corollary we obtain that the  relative quasiconvexity is equivalent
to the  dynamical quasiconvexity. This answers to a question posed
by D.~Osin \cite{Os06}.

In the second part of the paper we prove that a subgroup $H$ of a
finitely generated relatively hyperbolic group $G$ acts
cocompactly outside its limit set if and only if it is
(absolutely) quasiconvex and every its infinite intersection with
a parabolic subgroup of $G$ has finite index in the parabolic
subgroup.

We then obtain  a list of different  subgroup properties and
establish relations between them.

\end{abstract}
\maketitle
\markboth{V.~Gerasimov, L.~Potyagailo (\today)}{Quasiconvexity (\today)}

\def\tupl#1{\langle #1\rangle}
\def\b{\mathbf}
\def\s{\mathsf}
\def\d{\boldsymbol\delta}
\def\ti{\widetilde}
\def\act{\curvearrowright}
\def\l#1{\vskip3pt\noindent\bf Lemma \the\secNum.#1\sl.}
\newcount\secNum\relax\secNum=0
\def\Sec{\global\advance\secNum by 1\relax\the\secNum. }

\font\cFont=cmr9%
\def\note#1{\{{\cFont #1}\}}%
\let\g\mathfrak%
\section{Introduction}
\subsection{Results and history}
Let $G$ be a discrete group acting by homeomorphisms on a
compactum $X$ with the  {\it convergence property}, i.e. the
induced action on the space $\bold\Theta^3 X$ of subsets of $X$ of
cardinality $3$  is properly discontinuous. We say in this case
that $G$ acts {\it $3$-discontinuously} on $X$. Denote by $T$ the
limit set $\bold\Lambda G$ and by $A=\bold\Omega G$ the
discontinuity domain for the action $G\act X$. We have $X =
T{\sqcup}A$.

It is well-known that the action of a word-hyperbolic group on its
Gromov boundary has convergence property \cite{Gr87}, \cite{Tu94}.
However there are   convergence actions of groups that are not
Gromov hyperbolic: the actions of non-geometrically finite
Kleinian groups or those containing  parabolic subgroups of rank
at least 2; the actions of finitely generated groups on the space
of ends and on their Floyd boundaries \cite{Ka03}; the actions of
the groups of homeomorphisms of spheres, discontinuous outside a
zero-dimensional set \cite{GM87}.

An important class of groups  form  \it relatively hyperbolic
groups \rm(RHG for short). B. Bowditch \cite{Bo97} proposed a
construction of the ``boundary'' for such groups. This is a
compactum where the group acts with the convergence property. A.
Yaman proved that a group is RHG if it acts  on a metrisable
compactum $X$ {\it  geometrically finitely}, i.e. every point of
$X$ is either conical or bounded parabolic \cite{Ya04}.

 We call an action $G\act X$
{\it 2-cocompact} if the induced action on the space
$\bold\Theta^2X$ of distinct pairs is cocompact.

It follows from \cite{Ge09}, \cite{Tu98} that an action of  a
finitely generated group $G$ on a compactum $X$ is geometrically
finite if and only if it admits a 3-discontinuous and 2-cocompact
action on  $X.$  So we will further regard the existence of a
3-discontinuous and 2-cocompact action  as the definition of RHG.
Note that an advantage of this definition is that many results
known for finitely generated RHG remain valid for non-finitely
generated ones \cite{GP10}.

Recall that a subset $F$ of the Cayley graph of a group $G$ is
called   (absolute) \it quasiconvex\rm\ if every geodesic with the
endpoints in $F$ belong to a uniform neighborhood of $F$
\cite{Gr87}. Similarly a subset $F$ of the Cayley graph of $G$ is
called  {\it relatively quasiconvex}    if every
  geodesic   in the relative Cayley graph with endpoints in
$F$ belongs to a uniformly bounded neighborhood of $F$ in the
absolute Cayley graph (see Subsection \ref{relHull}).

 B.~Bowditch \cite{Bo99} characterized the
quasiconvex subgroups of Gromov hyperbolic groups in terms of
their action  on the Gromov boundary of the group.  He proved that
  a subgroup $H$ of a hyperbolic group $G$ is quasiconvex if and only if
   for any two disjoint closed subsets $K$ and $L$ of $T$ there are at most finitely many
distinct elements of $G$ such that the images of the limit set of
$H$ under them   intersect  both $L$ and $K$ (see Subsection 4.3).
He calls the latter  property  {\it dynamical quasiconvexity}.
 A natural question arises: can the dynamical
quasiconvexity for relatively hyperbolic groups be expressed  in
geometrical terms as it occurred for hyperbolic groups.
 D.~Osin  has
conjectured [Os06, Problem 5.3] that, for RHG, the relative
quasiconvexity is equivalent to  the Bowditch  dynamical
quasiconvexity. This conjecture follows from our first main
result.

We consider   other ``relative quasiconvexity'' properties.
Partially dynamical, partially geometrical. One of them is called
{\it visible quasiconvexity} and means that the set of points of
$A$ such that a given set $F\subset X$ has sufficiently big
diameter with respect to a shortcut metric (see 2.5)  based at
  points which must belong to a bounded neighborhood of $F$ with respect to the graph
  distance (see 4.3).

Generalizing the notion of relative quasiconvexity we call  a
subset $F\subset A$ {\it $\alpha$-relatively quasiconvex}    for
some
  distortion function $\alpha$ if every
$\alpha$-distorted path with endpoints in $F$ and  outside the
system of horospheres belongs to a bounded neighborhood of $F$
(see \ref{relHull}).

Our first main result shows that all these notions of the relative
quasiconvexity are equivalent.

\vskip5pt \bf Theorem A\sl. Let a finitely generated discrete
group $G$ act  3-discontinuously and 2-cocompactly on a compactum
$X$. The following properties of a subset $F$ of the discontinuity
domain of the action are equivalent:\hfil\penalty-10000
--- $F$ is relatively quasiconvex;\hfil\penalty-10000
--- $F$ is visibly quasiconvex;\hfil\penalty-10000
--- $F$ is relatively $\alpha$-quasiconvex where $\alpha$ is a quadratic polynomial with big enough coefficients.

Moreover, if $H$ is a subgroup of $G$ acting cofinitely on $F$
  then the visible quasiconvexity of $F$ is equivalent to the
dynamical quasiconvexity of $H$ with respect to the action
$G{\curvearrowright}X$\rm.\hfill $\blacklozenge$

\medskip

 Note that  in the first assertion of  Theorem A  we do not require that
  $F$ is acted upon by a subgroup of
 $G.$
 If in particular $F=H$  is a subset of the Cayley graph of $G$ then the
  second assertion  of the Theorem implies the following Corollary answering
 affirmatively   the above question of Osin.

\vskip5pt

 \bf Corollary\sl.\ A finitely generated subgroup $H$ of a
relatively hyperbolic group $G$ is relatively quasiconvex if and
only if it is dynamically quasiconvex\rm.

\vskip5pt

A graph $\Gamma$ acting upon by a group $G$ is called $G$-cofinite
if it has at most  finitely many $G$-non-equivalent edges; it is
called {\it fine} if for every $n{\in}\N$ and for every edge $e$
the set of simple loops in $\Gamma$ passant par $e$ of length $n$
is finite. As an application of the above methods we obtain a
generalization of the following result known for finitely
generated groups \cite{Ya04} to the case of infinitely generated
(in general uncountable) RHG.

\vskip5pt

\bf Proposition\ (\ref{genyam}).\sl\   Let $G$ be a group acting
2-cocompactly and 3-discontinuously on a compactum $T$. Then there
exists a hyperbolic, $G$-cofinite graph $\Gamma$ whose vertex
stabilizers are all finite except the vertices corresponding to
the parabolic points for the action $G\act T.$ Furthermore the
graph $\Gamma$ is fine.\rm

\vskip5pt

 The aim of the second
part of the paper is to relate different notions of  the
(absolute)  quasiconvexity of the subgroups of a relatively
hyperbolic group. We call a finitely generated subgroup $H$ of a
such group $G$ {\it weakly $\alpha$-quasiconvex} if $H$ acts
properly on $A$ (i.e. the point stabilizers are finite) and there
exists an orbit of $H$ for which every two points can be joined by
an $\alpha$-distorted path lying in a uniformly bounded
neighborhood of the orbit. This is a priori a partial case of a
more general  definition according to which $H$ is
$\alpha$-quasiconvex if {\bf every}  $\alpha$-distorted path in
$A$ connecting two points of an $H$-invariant and $H$-finite set
$E$ (i.e. $\vert E/H\vert <\infty$) is contained in a bounded
neighborhood of $E$.

The  following result describes the case when both these
conditions are equivalent to a stronger property to have cocompact
action outside the limit set.

\vskip5pt \bf Theorem B\sl.\ Let a finitely generated group $G$
act 3-discontinuously and 2-cocompactly on a compactum $X$. Let
  $\s{Par}$ be the set of the parabolic points for this action.
  Suppose that
$A{=}X\setminus T{\ne}\varnothing$ where $T{=}\b\Lambda G$ is the
limit set for the action. Then there exists a constant
$\lambda_0\in]1,+\infty[$ such that the following properties of a
subgroup $H$ of $G$ are equivalent:\hfil\penalty-10000 $\s a:$$H$
is weakly $\alpha$-quasiconvex for some
 distortion function $\alpha$ for which $\alpha(n){{\leqslant}}\lambda_0^n\ (n\in\N),$ and for
every $\g p{\in}\s{Par}$ the subgroup $H{\cap}\s{St}_G\g p$ is
either finite or has finite index in $\s{St}_G\g
p$;\hfil\penalty-10000 $\s b:$ the space $(X\setminus \b\Lambda
H)/H$ is compact;\hfil\penalty-10000 $\s c:$ for every distortion
function $\alpha$ bounded by $\lambda_0^n\ (n\in\N),$ every
$H$-invariant $H$-finite set $E{\subset}A$ is $\alpha$-quasiconvex
and  for every $\g p{\in}\s{Par}$ the subgroup $H{\cap}\s{St}_G\g
p$ is either finite or has finite index in $\s{St}_G\g
p$\rm.\hfill $\blacklozenge$

\vskip10pt

 The choice of the above constant $\lambda_0$ will be discussed
in \ref{Flmap}. In particular every subexponential function
satisfies our hypothesis. We also note that Theorem B shows that
the cocompactness outside the limit set is a stronger condition
than the usual quasiconvexity as it requires to preserve the
parabolic subgroups in the above sense. One of the applications of
the method used in the proof is the following.

\vskip5pt

\bf Proposition\ (\ref{finpres})\sl.\ Let $G$ be a group acting
3-discontinuously and 2-cocompactly on a compactum $X.$ Suppose
$H$ is a subgroup of $G$ acting cocompactly on $X\setminus
\bold\Lambda H$. If $G$ is finitely presented then $H$ also is.\rm

\vskip5pt

 Note that every maximal parabolic subgroup of a RHG acts
cocompactly outside its limit point on $X$. The above Proposition
is known   in the case when $H$ is maximal parabolic \cite{DG10}.
However it is easy to construct an  example  of  a quasiconvex
subgroup $H$ which cannot be parabolic for any geometrically
finite action of the ambiant group $G$  such that $H$  still
admits a cocompact action outside its limit set  (Example 1,
Subsection 9.2). We provide a direct proof of the Proposition in
this more general case.

 We note also that the cocompactness
 on the space $X\setminus \b\Lambda H$
differs from the cocompactness on the ``thinner" space $T\setminus
\b\Lambda H$ where $T{=}\b \Lambda G.$ There exist examples of
finitely generated discrete (Kleinian) subgroups of the isometry
group ${\rm Isom}\h3$ of the real hyperbolic space $\h3$ acting
non-geometrically finitely on $\h3$ and cocompactly outside their
limit sets on $\Bbb S^2$ (so called {\it totally degenerate
groups}). L.~Bers proved that they appear on the boundary of the
classical Teichm\"uller space of a closed surface \cite{Be70}.

 A subgroup $H$ of a group $G$ acting
3-discontinuously on $X$ is called {\it dynamically bounded} if
every infinite set $S$ of pairwise distinct elements of $G$ modulo
$H$ contains an infinite subset $S_0$ such that $T\setminus
\bigcup_{s\in S_0} s(\La H)$ has a non-empty interior (see
Proposition \ref{equivbound} for equivalent definitions). Our next
result is the following.

\vskip5pt

\bf Proposition\ (\ref{dcb})\sl.\  Let  $G$ act 3-discontinuously
on $X{=}T{\sqcup}A$. Suppose $H$ is a dynamically bounded subgroup
of $G$ acting cocompactly on $T\setminus \La H$. Then   $H$ acts
cocompactly on $\widetilde T\setminus \La H.$\rm

\vskip5pt

We note that we do not assume in the Proposition that the action
of $G$ on $X$ is 2-cocompact nor that $G$ is finitely generated.
We note also that the property opposite to the dynamical boundness
was studied by C.~McMullen \cite{McM96}. If a discrete  subgroup
$H$   is not dynamically bounded in   the full isometry group
$G{=}{\rm Isom}\h3$ then he says that  the limit set of $H$
contains {\it deep points} (i.e. {\it repeller}  points for
infinite sequences of elements of $G$  converging to a limit cross
and not satisfying the above definition). It turns out that the
limit set of a totally degenerate group contains uncountably many
deep points \cite[Corollary 3.15]{McM96}. It seems to be an
intriguing question to know whether  such an example of a finitely
generated subgroup of a relatively hyperbolic (or even
geometrically finite Kleinian) group could exist.

To summarize we obtain a list of  different properties of
subgroups of a relatively hyperbolic group. The following diagram
illustrates a natural order relation between them established in
the paper.

\vspace*{1cm}

%\begin{figure}[sh]
%\leavevmode{\epsfxsize=5in\epsfbox{figure20.eps}}
%\caption{Relations between the Theorems.}
%\end{figure}

\hspace*{0cm}{\epsfxsize=6.5in\epsfbox{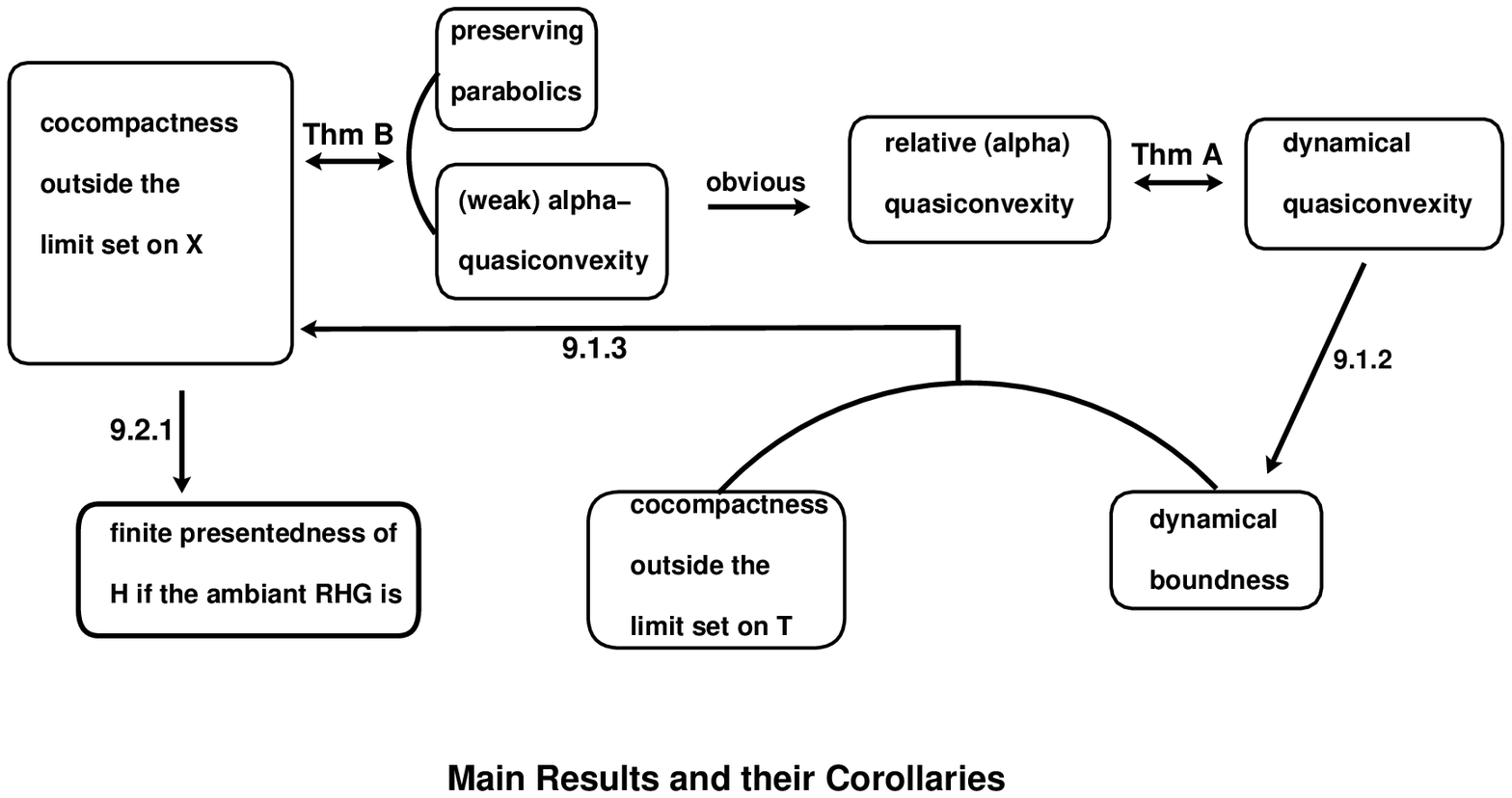}}

%\bigskip

%\centerline{Main Results and their Corollaries}

  \subsection {The structure of the paper} In Section 2 we generalize
   a useful lemma by A. Karlsson \cite{Ka03} about the Floyd length of ``far'' geodesics
    to the $\alpha$-distorted curves for some appropriate scalar function
$\alpha$ (Section 2). Using this lemma we obtain in  Section 3 a
uniform bound  for the size of the projections of subsets of $X.$

Generalizing the ideas of   \cite{GP09}
  we   prove in Section 4 that  the ($\alpha$-)convex
  hull of a closed set in $X$ is itself closed in $X$ (\ref{hullClosed}).
 As a corollary we
obtain that  the subgroups of relatively hyperbolic groups acting
cocompactly on $X$ outside  their limit sets are undistorted
(Corollary \ref{qIso}) and $\alpha$-quasiconvex (Proposition
\ref{hFinDiscr}). We then introduce  a notion of visible
quasiconvexity and prove that it is equivalent to the dynamical
quasiconvexity (\ref{vistodyn}).

In Section 5 we discuss the notion  of a general system of
horospheres. In particular we obtain  a uniform bound for the size
of the  projection of one horosphere onto another one
(\ref{horosph}). This result is used    in the sequel.

The proof of Theorem A is completed in Section 6. We first prove
that a lift of a geodesic path from the relative Cayley graph to
the absolute Cayley graph  is   $\alpha$-distorted for a quadratic
polynomial $\alpha$  (\ref{liftUndistorted}). We use it to prove
that a relatively quasiconvex subset is visibly quasiconvex
(\ref{relDyn}). These results imply Theorem A.

To illustrate the effectiveness of our methods we give in the
Section 7 simple independent proofs of some known results about
RHGs which use heavy techniques and require heavy references. A
new result obtained here is the above Proposition \ref{genyam}.

In Section 8 we prove all statements of Theorem B in the cyclic
order. The most difficult part is to prove the implication
  `$\s{a{\Rightarrow}b}$'. This is done by constructing a discrete
  analog of the Dirichlet fundamental polyhedron for a
  discrete group acting on $\hn.$ The main step is to prove that
  this set, denoted by $F_v\ (v{\in} A)$, is compact in
  $X\setminus\b\Lambda H$. The proof   is based on the
  methods developed in Section 4.

 In the last Section 9 we study subgroups of
convergence groups which admit the   dynamically boundness
property.   We prove here Propositions \ref{dcb} and \ref{finpres}
mentioned above. The dynamical boundness turns out to be  the
weakest subgroup property studied in the paper: all other
quasiconvexity properties imply it (see the table above).
 At the end of the Section we
provide some examples of dynamically bounded subgroups which are
not (relatively) quasiconvex and  not finitely presented. We
finish the paper by stating several questions which seem to be
open and intriguing.

\medskip
{\bf Acknowledgements.} During the work on this paper both authors
were partially supported by the ANR grant ${\rm BLAN}~07-2183619.$

The authors are thankful to Misha Kapovich and to Wenyuan Yang for
useful  discussions and suggestions.

\section{Karlsson functions for generalized quasigeodesics}

\subsection{Notations and definitions} We keep some notations and
terminology of \cite{GP09} and \cite{Ge10}. The canonical distance
function on the set $\Gamma^0$ of vertices of a graph $\Gamma$ is
denoted by $\s d$. By $\Gamma^1$ we denote the set of pairs of
vertices joined by edges.

For a subset $S$ of a metric space $(M;\d)$ and a nonnegative
number $r$ we consider the\hfil\penalty-10000
$r$\it-neighborhood\rm\ $\s N_r^{\d}S{\leftrightharpoons}\{\g
p{\in}M:\d(S,\g p){\leqslant}r\}$. For a set $S$ of vertices of a
graph we sometimes write $\s N_rS$ instead of $\s N_r^{\s d}S$.

For a path $\gamma:I\to\Gamma^0$ in a graph $\Gamma$ we call $I$
its \it domain \rm$\s{Dom}\gamma$ and the set $\gamma(I)$ its \it
image \rm$\s{Im}\gamma$. The diameter of $\s{Dom}\gamma$ is the
\it length \rm of $\gamma$. If $|I|{<}\infty$ we define
$\partial\gamma{\leftrightharpoons}\gamma(\partial I)$. We extend
naturally the meaning of $\partial\gamma$ over the half-infinite
and bi-infinite paths in the case when $\Gamma^0$ is a discrete
subset of a Hausdorff topological space $X$ and the corresponding
infinite branches of $\gamma$ converge to points of $X$.

By $\s{length}_{\d}\gamma$ we denote the length of a path $\gamma$
with respect to a path-metric $\d$. \vskip3pt Considering a
function $f$ defined on a subset of $\Bbb Z$ we sometimes write
$f_n$ instead of $f(n)$.

By $|S|$ we denote the cardinality of a set $S$. By
$\bold\Theta^nS$ we denote the set of all subsets of $S$ of
cardinality $n$. When $S$ is a topological space then
$\bold\Theta^nS$ is considered with the induced topology. By $\s
S^nS$ we denote the set of ``generalized unordered $n$-tuples'':
formally this is the quotient of the Cartezian power $S^n$ by the
action of the permutation group.

If $S$ is acted upon by a group $G$ it is called $G$-set. A subset
 $M$ of $G$-set $S$ is called $G$\it-finite \rm if $M$ meets
finitely many $G$-orbits. In this case the image of $M$ in $S/G$
is finite.

Recall that a limit point $p\in\bold\Lambda G$ for the convergence
action of $G$ on a compactum $S$ is called {\it parabolic} if it
is the unique limit point for the action of its stabilizer
$\s{St}_Gp{=} \{g\in G\ :\ gp{=}p\}$   on $S.$ A parabolic limit
point $p\in \bold\Lambda G $ is called
  {\it bounded parabolic} if $S\setminus\{p\}/\s{St}_Gp$ is
compact.

\subsection{Distorted paths} A nondecreasing function $\alpha:\Bbb N\to\Bbb R_{>0}$ such that
$\forall n\ \alpha_n{\geqslant}n$ is called a \it distortion
function\rm. Thus, the minimal distortion function is the function
$\s{id}:n\mapsto n$.

Let $\alpha$ be a distortion function. A path
$\gamma:I\to\Gamma^0$ in a graph $\Gamma$ is said to be
$\alpha$\it-distorted \rm if
$\s{diam}J{\leqslant}\alpha(\s{diam}\gamma(\partial J))$ for every
\bf finite \rm interval $J{\subset}I$.

If $\alpha$ has one of the following forms (1) $n\mapsto n$,
(2)$n\mapsto Cn$, (3) $n\mapsto Cn{+}D$, the notion
`$\alpha$-distorted' means respectively `geodesic', `Lipschitz',
`large-scale Lipschitz' (`quasigeodesic'). The case when $\alpha$
is a quadratic polynomial will be of our particular interest.
%\vskip3pt Note that if $\gamma$ is a non-bi-infinite path and
%$\alpha_0<\s{length_d}\gamma$ then $\partial\gamma$ is a proper
%pair. If $\gamma$ is bi-infinite and $\partial\gamma$ is a single
%point $\g p$, we call $\gamma$ an $\alpha$\it-loop \rm at $\g p$.
\subsection{Scaling of the graph metric}\label{scale}
Recall the notions related to the Floyd metrics.

A function $f:\Bbb N\to\Bbb R$ is said to be a (Floyd) \it scaling function \rm
if $\sum_{n\geqslant0}f_n<\infty$ and there exists a positive $\lambda$ such that $1\geqslant f_{n+1}/f_n\geqslant\lambda$ for all $n{\in}\Bbb N$.
The supremum of such numbers $\lambda$ is called the \it decay rate \rm of $f$.

Let $f$ be a scaling function and let $\Gamma$ be a connected graph. For each vertex $v{\in}\Gamma^0$ we define on $\Gamma^0$ a new
metric $\d_{v,f}$ as the maximal among the metrics $\varrho$ on $\Gamma^0$ such that $\varrho(x,y)\leqslant f(\s d({v,\{x,y\}}))$
for each $\{x,y\}{\in}\Gamma^1$.
We say that $\d_{v,f}$ is the \it Floyd metric \rm(with respect to the scaling function $f$) \it based \rm at $v$.

When $f$ is fixed we write $\d_v$ instead of $\d_{v,f}$. When $v$ is also fixed we write $\d$ instead of $\d_v$.

One verifies that $\d_u/\d_v\geqslant\lambda^{\s d(u,v)}$ for
$u,v{\in}\Gamma^0$. Thus the Cauchy completion $\overline\Gamma_f$
of $\Gamma^0$ with respect to $\d_{v,f}$ does not depend on $v$.
The \it Floyd boundary \rm is the space
$\partial_f\Gamma{\leftrightharpoons}\overline\Gamma_f\setminus
\Gamma^0$. Every $\s d$-isometry of $\Gamma$ extends to a
homeomorphism $\overline\Gamma_f\to\overline\Gamma_f$. The Floyd
metrics extend continuously onto the Floyd completion
$\overline\Gamma_f$.
\subsection{Karlsson functions}\label{Ka}
Let $f$ be a Floyd function and let $\alpha$ be a distortion
function.\hfil\penalty-10000 A non-increasing function $\s K:\Bbb
R_{>0}\to\Bbb N$ is called \it Karlsson function \rm for the pair
$(f,\alpha)$ if
\begin{equation}\label{KaFn}
\s d(v,\s{Im\gamma)\leqslant\s K(length}_{\d_{v,f}}\gamma)
\end{equation} for each $\alpha$-distorted path $\gamma$ in a connected graph with a vertex $v$.
A pair $(f,\alpha)$ where $f$ is a scaling function and $\alpha$
is a distortion function is said to be \it appropriate \rm if it
possesses a Karlsson function.  It is proved in \cite{Ka03} that
every pair of the form $(f,\s{id})$ is appropriate. A similar
agrument can be applied to show that $(f,\alpha)$ is appropriate
for $\alpha:n\mapsto Cn{+}D$.

We need one more class of appropriate pairs. Actually, all pairs considered in this article belong to this class.

\begin{prop}\label{gKaFn}
If $\sum_{n\geqslant0}\alpha_{2n+1}f_n<\infty$ then the pair $(f,\alpha)$ is appropriate.
\end{prop}
\it Proof\rm. Let $v$ be the reference point of some graph
$\Gamma$. Denote $|x|{\leftrightharpoons}\s d(v,x)$.

Let $\gamma:I\to \Gamma^0$ be an $\alpha$-distorted path. We can
assume that $\s
d(v,\s{Im}\gamma){=}|\gamma(0)|{\leftrightharpoons}r$. It suffices
to prove that $\s{length}_{\d}(\gamma|_{I{\cap}\Bbb N})$ is small
enough whenever $r$ is big enough. So we can assume that $I$ is an
initial segment of $\Bbb N$.

By induction we define a strictly increasing sequence $x_s{\in}I$
for $s{\geqslant}r$ such that $|\gamma(x_s)|{=}s$ and
$\gamma([x_s,x_{s+1}]){\cap}\s N_{s-1}^{\s d}v{=}\varnothing$.
Indeed let $x_r{\leftrightharpoons}0$. If $x_s$ is already defined
and different from $\s{max}I$, put $x_{s+1}\leftrightharpoons
1{+}\s{max}\{x{\in}I:x{\geqslant}x_s$ and $\s
d(v,\gamma(x)){=}s\}$. Now the interval $I$ has subdivided into
the segments $I_s{\leftrightharpoons}[x_s,x_{s+1}]$. By the
$\triangle$-inequality we have
$\s{length}_{d}\gamma|_{I_s}{=}\s{diam_d}I_s {=}
x_{s+1}{-}x_s{\leqslant}\alpha_{2s{+}1},$ hence
$\s{length}_{\d}\gamma|_{I_s}\leqslant f(d(v, I_s))\cdot
\s{length}_{d}\gamma|_{I_s} \leqslant \alpha_{2s{+}1}f_s$. Thus
$\s{length}_{\d}\gamma\leqslant\sum_{s=r}^{k-1}\alpha_{2s+1}f_s+\alpha_{2k}f_k$
and $k < +\infty$   only if $I$ is finite and
$\vert\gamma(x_k)\vert{=}\vert\gamma(\s{max}I)\vert$. In any case
we have
\begin{equation}\label{resser}
\s{length}_{\d}\gamma\leqslant\sum_{s=r}^{\infty}\alpha_{2s+1}f_s\
\s{where}\ r{=}\s d(v,\s{Im}\gamma).\end{equation}

 Thus the function
$\varepsilon\mapsto\s{min}\{r:\sum_{s=r}^\infty\alpha_{2s{+}1}f_s\leqslant\varepsilon/2\}$
is a Karlsson function for $(f,\alpha)$. Indeed if not then
$r{=}\s {d}(v,\s{Im}\gamma)>\s {K(length}_{\d}\gamma){=}r_0$ and
$\sum_{s=r}^{\infty}\alpha_{2s+1}f_s\leqslant\sum_{s=r_0}^{\infty}
\alpha_{2s+1}f_s<\s{length}_{\d}\gamma$ contradicting
\ref{resser}.\qed

\medskip

It follows immediately from \ref{gKaFn} that every
$\alpha$-distorted ray converges to a point at the Floyd boundary.
So, an $\alpha$-distorted ray extends to a continuous map
$\overline{\Bbb N}{\leftrightharpoons}\Bbb
N{\cup}\infty\to\overline\Gamma_f$. Similarly, any
$\alpha$-distorted line extends to a continuous map
$\overline{\Bbb Z}{\leftrightharpoons}\Bbb
Z{\cup}\{{\pm}\infty\}\to\overline\Gamma_f$ (see [GP09,
Proposition 2.5] for the case of affine distortion).

Now $\partial\gamma{\leftrightharpoons}\gamma(\partial I)$ is
well-defined for finite, half-infinite and bi-infinite
$\alpha$-distorted paths. It is a subset of
$\overline\Gamma_\lambda$.
\subsection{Floyd map}\label{Flmap}
 From now on we fix a compactum which we denote for the sake of convenience by
 $\widetilde T$.  We also fix a 3-discontinuous action of a
 discrete group $G$ on $\widetilde T$. If the opposite is not stated we will
also suppose that the action is 2-cocompact. We have $\widetilde T
{=} T{\sqcup}A$ where $T{=}\bold\Lambda G$ is the limit set and
$A{=}\bold\Omega\widetilde T$ is the discontinuity set for the
action $G{\curvearrowright}\widetilde T$.     Up to adding a
discrete $G$-orbit to the space $T$ we can always assume that $A$
is a non-empty, discrete and  $G$-finite set (i.e. $\vert A/G\vert
<\infty$), and the compactum $\ti T$ contains at least 3 points
(i.e. $\bold\Theta^3 \ti T{\not=}\emptyset$).

Let $\Gamma^1$ be a $G$-finite subset of $\bold\Theta^2A$ such
that the graph $\Gamma$ with $\Gamma^0{=}A$ is connected. Such
subset exists if and only if $G$ is finitely generated (see e.g.
\cite{GP09} or \cite{GP10}  where $A$ is the vertex set of the
Cayley graph of $G$ or a $G$-orbit of  an entourage of $T$).

\medskip

 \bf Convention 1\rm.    Since now on by default we assume that
$G$ is a finitely generated group acting on a locally finite,
$G$-finite and  connected graph $\Gamma$ such that $\Gamma^0{=}A.$
 We will
always assume that $\vert \ti T\vert > 2.$
\medskip

 It is proved
in \cite{Ge10} that there exists an exponential scaling function
$f_0(n){=}\mu_0^n\ \hfil\penalty-10000 (\mu_0{\in}]0,1[,\
n{\in\N}),$
 and a metric $\varrho$ on $\widetilde T$
determining the topology of $\widetilde T$ such that
$\d_{v,f_0}{\geqslant}\varrho$ on $A$ where $\d_{v,f_0}$ is the
Floyd metric on $\overline\Gamma_{f_0}$ at a point $v\in A.$ Thus
the inclusion map $A\hookrightarrow\widetilde T$ extends
continuously to the map $\overline\Gamma_{f_0}\to\widetilde T$
called \it Floyd map\rm.

 The Floyd map  induces a set of   shortcut
metrics $\overline\d_v $ on $\widetilde T \ (v{\in} A)$,  where
every $\overline\d_v$ is the maximal   among all  metrics $\rho$
on $\widetilde T$ \cite{GP09}.

We denote by $\lambda_0\in ]1, +\infty[$ the maximal  constant for
which the distortion function  $\alpha_n{=}\lambda_0^n\ (n\in\N)$
is appropriate for the above  Floyd function $f_0.$

\medskip

\vskip3pt \bf Convention 2\rm. We will always consider a Floyd
function $f$ satisfying $f(n){\geqslant}f_0(n)\ (n\in\N)$, so the
Floyd map also exists for $f.$ For a fixed Floyd function
  $f$ we will always choose an appropriate
distortion function $\alpha$ ($\alpha_n{\leqslant}\lambda_0^n\
(n\in \N)).$

 For every appropriate pair
$(f,\alpha)$ we fix a Karlsson function denoted by
 $\s K_{f,\alpha}$. We also write $\s K_\alpha$ instead of $\s K_{f,\alpha}$ and $\s
K$ instead of $\s{K_{id}}$.

\section{Projections}
\subsection{Boundary equivalence}\label{bEqu}
For a set $E{\subset}A$ define $\partial
E{\leftrightharpoons}T{\cap}\overline E$. This ``boundary'' is
nonempty if and only if $E$ is infinite. Since $A$ is a discrete
open subset of a compactum, for any neighborhood $N$ of $\partial
E$ in $\widetilde T$ the set $E\setminus N$ is finite. In
particular $\overline E{=}E{\cup}\partial E$. Thus, for $a{\in}A$
and $\varepsilon{>}0$  the number \setcounter{equation}0
\begin{equation}\label{CEa}\s
C_{E,a}(\varepsilon){\leftrightharpoons}\s{min}\{r:E\setminus \s
N_r^{\s d}a\subset\s N_\varepsilon^{\overline\d_a}\partial
E\}\end{equation}

\noindent is finite, where $\s N_r^{\s d}$ and $\s
N_\varepsilon^{\overline\d_a}$ are $r$ and
$\varepsilon$-neighborhoods with respect to the metrics $d$ on $A$
and $\overline\d_a$ on $\widetilde T$ respectively.

\medskip

 {\bf Definition.} Two sets $E,F{\subset}A$ are said to be
$\partial$\it-equivalent \rm(notation $E{\sim}_\partial F$) if
$\partial E{=}\partial F$.

\medskip
\setcounter{prop}1
\begin{prop}\label{bEquN} $E{\sim}_\partial\s N_r^{\s d}E$ for every $E{\subset}A$, $r{\in}\Bbb N$.\end{prop}
\it Proof\rm. It suffices to prove the statement for $r{=}1$. The result follows from the fact that the metric $\overline\d_a$
determines the topology of $\widetilde T$ and that the $\overline\d_a$-length of an edge $e$ tends to zero while
$\s d(a,e)\to\infty$.\qed

\subsection{Projections of subsets of $A$}
For a vertex $a{\in}A$ define the \it projection set
\rm$\s{Pr}_Ea{\leftrightharpoons}\{v{\in}E:\s d(a,v){=}\s
d(a,E)\}$. For a nonempty set $B{\subset}A$ define
$\s{Pr}_EB{\leftrightharpoons}{\cup}\{\s{Pr}_Eb:b{\in}B\}$.

\begin{prop}\label{prOfA}
If  $\partial E{\cap}\partial B{=}\varnothing$ $(E,B{\subset}A)$ then $\s{Pr}_EB$ is finite.
\end{prop}
\it Proof\rm. Suppose that $\partial E{\cap}\partial
B{=}\varnothing$ for $E,B{\subset}A$. We can assume that $E$ is
infinite and hence $\partial E{\ne}\varnothing$. Since $\partial
E{\cap}\overline B{=}\varnothing$ the number

\begin{equation}\label{numbrho}\rho{\leftrightharpoons}\rho(E,B){\leftrightharpoons}
\s{sup}\{\overline\d_v(\partial E,\overline
B):v{\in}E\}\end{equation} is positive. Let
$0{<}\delta{<}\varepsilon{<}\rho$ and let $a{\in}E$ be such that
$\overline\d_a(\partial E,\overline B){>}\varepsilon$. We will
show that $\s{Pr}_EB$ is within a bounded distance from $a$.

Denote $r{\leftrightharpoons}\s C_{E,a}(\varepsilon{-}\delta)$. If
$b{\in}B$, $v{\in}\s{Pr}_Eb$ then either $v{\in}\s N_r^{\s d}a$ or
$v{\in}\s N_{\varepsilon-\delta}^{\overline\d_a} \partial E.$ In
the latter case we have $\overline\d_a(b,v)
{\geqslant}\overline\d_a(\partial E,\overline B)-\overline\d_a(v,
\partial E){>}\varepsilon -\varepsilon +\delta{=}\delta$. Thus
 for a geodesic segment $\gamma$ between $b$ and $v$ by \ref{KaFn}
we have $\s d(a,\s{Im}\gamma){\leqslant}s{\leftrightharpoons}\s
K(\delta)$. Therefore for $c{\in}\s{Im}\gamma$ such that $\s
d(a,c){=}\s d(a,\s{Im}\gamma)$ we obtain  $\s d(b,c){+}\s
d(c,v){=}\s d(b,v){\leqslant}\s d(b,a){\leqslant} \s d(b,c){+}\s
d(c,a)$. Thus $\s d(c,v){\leqslant}\s d(c,a)$ and
\hfil\penalty-10000 $\s d(a,v){\leqslant}\s d(a,c){+}\s
d(c,v){\leqslant}2s$. It yields
\begin{equation}\label{diampr}d(a,v){\leqslant}\s 2\cdot \max\{r,
\s K(\delta)\}.\end{equation}\qed
\subsection{Projection of the subsets of $\widetilde T$}\label{prExtent}
For a set $F{\subset}\widetilde T$ denote by $\s{Loc}_{\widetilde
T}F$ the set of all neighborhoods of $F$ in $\widetilde T$.

Let $E{\subset}A$. A $\widetilde T$-neighborhood $P$ of a point
$\g p{\in}T\setminus \partial E$ is called $E$\it-stable \rm if
$\s{Pr}_E(P{\cap}A){=}\s{Pr}_E(Q{\cap}A)$ for every
$Q{\in}\s{Loc}_{\widetilde T}\g p$ such that $Q{\subset}P$. By
\ref{prOfA} every point $\g p{\in}T\setminus \partial E$ possesses
a $E$-stable neighborhood since otherwise we would have a strictly
decreasing infinite sequence of sets of the form
$\s{Pr}_E(P{\cap}A)$, $P{\in}\s{Loc}_{\widetilde T}\g p$. If $P,Q$
are $E$-stable neighborhood of $\g p$ then $P{\cap}Q$ is also an
$E$-stable neighborhood and
$\s{Pr}_E(P{\cap}A){=}\s{Pr}_E(P{\cap}Q{\cap}A){=}\s{Pr}_E(Q{\cap}A)$.

Now we can extend the projection map over $\widetilde T\setminus
\partial E$: the projection $\s{Pr}_E\g p$ of a point $\g
p{\in}T\setminus \overline E$ is the projection of any its
$E$-stable neighborhood.

%\begin{prop}\label{eStable}
%Fvery open in $\widetilde T$ subset of $T\setminus \overline E$ is a union of $E$-stable sets.
%\end{prop}
\vskip3pt We need a uniform estimate for the size of the
projection. To this end we put \setcounter{equation}0
\begin{equation}\label{CE}
 C_E(\varepsilon){\leftrightharpoons}\s{sup\{C}_{E,a}(\varepsilon):a{\in}E\}.
 \end{equation}

 Let us call an \bf infinite \rm set $E{\subset}A$
\it weakly homogeneous \rm
 if $C_E(\varepsilon)<\infty$
for every $\varepsilon{>}0$.

 The following example is motivating.
Let $H$ be an infinite subgroup of $G$ and let $E$ be an
$H$-finite subset of $A$. Since $\s C_{E,a}(\varepsilon){=}\s
C_{gE,ga}(\varepsilon)\ (g{\in}G)$ the set
$\{C_{E,a}(\varepsilon):a{\in}E\}$ is finite for every
$\varepsilon>0$. Hence $E$ is weakly homogeneous.

We call $\s C_E$ the \it convergence function \rm for $E{\subset}A$. Its role is similar to that ot Karlsson functions.

Assuming that the constants $\varepsilon $ and $\delta$ from the
proof of \ref{prOfA} satisfy $\delta{\geqslant}\rho/4$ and
$\varepsilon -\delta{=}\rho/4$ we have
\setcounter{prop}1%
\begin{prop}\label{prUniSize}
For a weakly homogeneous set $E$
 the $\s d$-diameter of $\s{Pr}_EB$ depends only on
the number $\rho{=}\rho(E,B)$ of \ref{numbrho} and the function
$\s C_E$. More precisely, \setcounter{equation}2%
\begin{equation}\label{projdist}
\s{diam_dPr}_EB\leqslant2{\cdot}\s{max\{C}_E(\rho/4),\s
K(\rho/4)\}.
\end{equation}\end{prop}\qed

 We extend the distance function $\s
d$ over $\widetilde T^2$ by setting $\s
d(\g{p,q}){\leftrightharpoons}\infty$ for $\g{q{\ne}p}$ and
$\g{p}{\in}T$. So, for $F{\subset}\widetilde T$ we have
 $\s N_r^{\s d}F{=}\partial F\cup\s N_r^{\s
d}(F{\cap}A)$.

\section{Hulls and convexity}
\subsection{$\alpha$-quasiconvexity and $\alpha$-hull}\label{aHull}
Let $\alpha$ be a distortion function. The $\alpha$\it-hull \rm of
set $F{\subset}\widetilde T$ is\hfil\penalty-10000 $\s H_\alpha
F{\leftrightharpoons}{\bigcup}\{\s{Im}\gamma:\gamma$ is an
$\alpha$-distorted path in $A$ and
$\partial\gamma{\subset}F\}$.\hfil\penalty-10000 A set
$F{\subset}\widetilde T$ is said to be $\alpha$\it-quasiconvex \rm
if $\s H_\alpha F\subset\s N_rF$ for some $r{<}\infty$. In the
case when $\alpha{=}\s{id}$ ``$\alpha$-quasiconvex'' means
``quasiconvex''.

In the sequel we will always assume that $\alpha$ satisfies the
hypothesis of Proposition \ref{gKaFn}. Since $\sum_{n\geqslant0}
f_n<+\infty$ the function $\alpha{+}1$ also satisfies it. On the
other hand, $\s N_r\s H_\alpha E{\subset}\s H_{\alpha+2r}E$. This
implies that \setcounter{equation}0\relax
\begin{equation}\label{hUnion}A\subset
{\cup}_{r\geqslant0}\s H_{\alpha+r}E \end{equation} for every
$E\subset A.$
\setcounter{prop}1
\begin{prop}\label{hullLoc}
For every $\varepsilon{>}0$ there exists a number
$s{=}s(\varepsilon,\alpha)$ such that $\s H_\alpha F{\subset}\s
N_\varepsilon^{\overline\d_a}F{\cup}\s N_s^{\s d}a$ for every
$F{\subset}\widetilde T$ and $a{\in}A$.
\end{prop}
\it Proof\rm. It is similar to that of \cite[Main Lemma]{GP09}.

Define $r{\leftrightharpoons}\s K_\alpha(\varepsilon)$, $s{\leftrightharpoons}r{+}{1\over2}\alpha(2r)$.

Let $v{\in}\s H_\alpha F\setminus \s
N_\varepsilon^{\overline\d_a}F$ and let $\gamma:I\to A$ be an
$\alpha$-distorted path with $\gamma(0){=}v$ and
$\partial\gamma{\subset}F$. Denote
$\gamma_+{\leftrightharpoons}\gamma|_{I\cap\Bbb N}$,
$\gamma_-{\leftrightharpoons}\gamma|_{I\cap(-\Bbb N)}$. Since
$\s{length}_{\d_a}\gamma_\pm\geqslant\d_a(F,v)\geqslant\overline\d_a(F,v)>\varepsilon$
 by \ref{KaFn} we have $\s d(a,\s{Im}\gamma_\pm){\leqslant}r$. Let
$J{\ni}0$ be a subsegment of $I$ with $\gamma(\partial
J){\subset}\s N_r^{\s d}a$. So $\s{diam_d}\gamma(\partial
J){\leqslant}2r$ and $\s{length_d}\gamma|_J{\leqslant}\alpha(2r)$.
Hence $d(v,a){\leqslant}d(a,\s{Im}\gamma_\pm) + {1\over
2}\cdot\s{length_d}\gamma|_J{\leqslant}r+{1\over
2}\alpha(2r){=}s$. So $v{\in}\s N_s^{\s d}a$. \qed
\begin{prop}\label{hullClosed}
$E\sim_\partial\s H_\alpha\overline E$ for every $E{\subset}A$.

\end{prop}
\it Proof\rm. Suppose by contradiction that $\g p{\in}\partial{\s
H_\alpha\overline E}\setminus \partial E$. If
$0{<}\varepsilon{<}\overline\d_a(\partial E,\g p)$ for $a{\in}A$
then, by \ref{hullLoc}, $\s H_\alpha E$ is contained in the closed
set $\s N_\varepsilon^{\overline\d_a}E{\cup}\s N_s^{\s d}a$ that
does not contain $\g p$. A contradiction.\qed
\subsection{Subgroups acting cocompactly outside its limit set}
%\begin{prop}\label{boundedClosed}
%If $F{\in}\s{Closed}\widetilde T$ and $F{\subset}E{\subset}\s N_r^{\s d}F$ for some $r{<}\infty$ then $E{\in}\s{Closed}\widetilde T$.
%\end{prop}
%\it Proof\rm.
%Let $\g p{\in}\overline E\setminus E$. Since the points of $A$ are isolated, we have $\g p{\in}T\setminus F$.
%Let $P$ be an $E$-stable (see \ref{prExtent}) neighborhood of $P$ and let $Q{\subset}P$ be a neighborhood
%that does not meet the finite set $\s N_r^{\s d}\s{Pr}_F\g p{=}\s N_r^{\s d}\s{Pr}_FQ$. We have $\s d(Q,F){=}\s d(Q,\s{Pr}_FQ){>}r$,
%$Q{\cap}E{=}\varnothing$.
%A contradiction.\qed
We provide below several properties of a subgroup $H$ of $G$
acting cocompactly on the complement $\widetilde T\setminus
\bold\Lambda H$ of its limit set $\Lambda (H).$ In particular the
group $H$ can be a parabolic subgroup of $G$ for the action on
$\widetilde T.$ However there are a lot of examples of subgroups
satisfying this property and which are essentially non-parabolic
(see Example 1 in Subsection \ref{finpres}). In the following
Proposition we use the projection map $\s{Pr}_E$ on a subset
$E\subset A$ introduced in 3.2.
\begin{prop}\label{qIso}
Let $H$ be a subgroup of $G$ acting cocompactly on $\widetilde
T\setminus \bold\Lambda H$ and let $E$ be a nonempty $H$-finite
$H$-invariant set. Then the multivalued map $\s{Pr}_E$ is
quasi-isometric i.e. there exists a constant $C$ such that
$\s{diam(Pr}_Ee)\leqslant C$ for every edge $e{\in}\Gamma^1$.
\end{prop}
\it Proof\rm. By 3-discontinuity, $\partial E{=}\bold\Lambda H$.
Let $K$ be compact set such that $HK{=}\widetilde T\setminus
\bold\Lambda H$. By \ref{bEquN} the set $K_1{\leftrightharpoons}\s
N_1^{\s d}K$ is closed and disjoint from $\partial E$. For every
edge $e{\in}\Gamma^1$ there exists $h{\in}H$ such that
$he{\subset}K_1$.   Then
$\s{diam(Pr}_Ee){=}\s{diam(Pr}_E(he)){\leqslant}C{\leftrightharpoons}\s{diam(Pr}_EK_1)$.

Since $E$ is weakly homogeneous by \ref{projdist} we have
$C<\infty.$ \qed

%\note{It should be interesting to generalize this from $H$-finite sets to the sets quasiisometric to $H$-finite sets}
\vskip3pt \bf Corollary\sl.\ If $H$ acts cocompactly on
$\widetilde T\setminus \bold\Lambda H$ then $H$ is undistorted in
$G$\rm.\qed

\medskip

\label{cocompNonLim}\label{hFinDiscr} If, for a subgroup $H{<}G$,
the space $(\widetilde T\setminus \bold\Lambda H)/H$ is compact
then every closed $H$-invariant set $E{\subset}A$ is $H$-finite.
Indeed, its image in $(\widetilde T\setminus \bold\Lambda H)/H$ is
a closed discrete subset of a compact space.

In \cite[Lemma 3.3]{GP09} we make use of this observation for
parabolic subgroups. The following Proposition show that the
stronger quasiconvexity property is true for such subgroups.

\begin{prop}\label{subgrConvex}For every subgroup $H{<}G$ acting cocompactly on $\widetilde T\setminus \bold\Lambda H$ and
every $H$-finite $H$-invariant set $E$ the set $\s
H_\alpha\overline E$ is $H$-finite and $\overline E$ is
$\alpha$-quasiconvex. In particular, $H$ is an
$\alpha$-quasiconvex subgroup for any appropriate distortion
function $\alpha$.\end{prop} \it Proof\rm. By \ref{hullClosed} the
set $\partial E{\cup}\s H_\alpha\overline E$ is closed. Hence the
$H$-invariant set $\s H_\alpha\overline E$ is closed in
$\widetilde T\setminus \partial E$. By the above observation it is
$H$-finite. Thus   $\s H_\alpha \overline E{\subset}\s N_r^{\s
d}E$ for some $r>0.$ So $\overline E$ is $\alpha$-quasiconvex.

 In the case when $E$ is a single $H$-orbit this means
the last statement of \ref{subgrConvex}.\qed

\subsection{%$\varepsilon$-hull
Dynamical and visible quasiconvexity}\label{eConvex} For a set
$F{\subset}\widetilde T$ define its $\varepsilon$\it-hull \rm as
$$\s V_\varepsilon
F{\leftrightharpoons}\{a{\in}A:\s{diam}_{\overline\d_a}F{\geqslant}\varepsilon\}.$$
Note that $\s V_\varepsilon F{=}\s V_\varepsilon\overline F$ for
every $F{\subset}\widetilde T$.

\medskip

\bf Definition\rm. A set $F{\subset}\widetilde T$ is said to be
\it visibly quasiconvex \rm if for every $\varepsilon{>}0$ there
exists $r{=}r(\varepsilon){<}\infty$ such that $\s V_\varepsilon
F\subset\s N_r^{\s d}F$. We call the function $\s
Q_F:\varepsilon\mapsto r(\varepsilon)$
  \it visible quasiconvexity function\rm.

   \medskip

\begin{prop}\label{eHull}
$\partial\s V_\varepsilon\overline E\subset\partial E$ for every $E{\subset}A$ and $\varepsilon{>}0$.
\end{prop}
\it Proof\rm. Let $v{\in}\s V_\varepsilon\overline E{=}\s
V_\varepsilon E$. Let $\gamma:I\to A$ is a $\s d$-geodesic with
$\partial\gamma{\subset} E$ and
$\s{diam}_{\overline\d_v}\partial\gamma\geqslant\varepsilon$. We
have  $\s d(v,\s{H_{id}}E)\leqslant\s d(v,\s{Im}\gamma)\leqslant
r{\leftrightharpoons}\s K(\varepsilon)$. Thus $\s
V_\varepsilon\overline E\subset\s N_r\s{H_{id}}E$. Since $\s
N_r\s{H_{id}}E\overset{\text{by \ref{bEquN}}}{\sim_\partial}
\s{H_{id}}E\overset{\text{by \ref{hullClosed}}}{\sim_\partial}E$
it yields $\partial \s V_{\varepsilon}\overline E\subset\partial
E.$\qed

\medskip

\bf Corollary\sl.\ Every quasiconvex set is visibly
quasiconvex.\hfil\penalty-10000 \it Proof. \rm It follows
immediately from the inclusions $\s V_\varepsilon\overline
E\subset\s N_r\s{H_{id}}E$ and $\s{H_{id}}E\subset \s N_{r_0}(E)$
where $r{=}\s K(\varepsilon)$ and $r_0$ is the quasiconvexity
constant of $E$.\qed
\medskip

We now recall  few facts needed in the sequel. Let $X$ be a
compactum. A neighborhood of the diagonal ${\bf \Delta}^2X$ of the
space $X^2$   is called \it entourage\rm. The set of all
entourages on $X$ is denoted by $\s{Ent}X$. For an entourage
$\bold u$ a subset $S$ of $X$ is called {\it $\bold u$-small} if
$S^2\subset \bold u$. The set of all $\bold u$-small sets is
denoted by $\s{Small}(\bold u)$ (see \cite{Ge09} and \cite{GP10}
for more details).

\medskip

\bf Definition\rm\  \cite{Bo99}.   A subgroup $H$ of a discrete
group $G$ acting 3-discontinuously on a compactum $X$ is said to
be \it dynamically quasiconvex \rm if for every entourage $\bold
u$ of $X$ the set $G_{\bold u}{=}\{g{\in}G:g(\bold\Lambda
H){\notin}\s{Small}(\bold u)\}$ is $H$-finite with respect to the
$H$-action from the \bf right \rm.

\medskip

\bf Remarks\rm.  a)  The above definition coincides with the
notion of the dynamical quasiconvexity proposed in \cite{Bo99}.
The latter one states   that the set of the left cosets
\begin{equation}\label{dyncon}\{gH\ :\ gS\cap
L{\not=}\emptyset\ {\rm and}\ gS\cap K{\not=}\emptyset\}
\end{equation}

\noindent is finite, whenever $K$ and $L$ are disjoint closed
subsets of $X$ and   $S{=}{\bold \Lambda} H.$  Note first that one
can consider here only the entourages of a special form $\bold
u_{P,Q}{=}\s S^2T\setminus P{\times}Q$ where $P$ and $Q$ are
disjoint closed sets  is not a restriction since  the set of
entourages of this form generates the filter $\s{Ent}T$ of
entourages.

In order to see that (\ref{dyncon}) is equivalent to the
definition above suppose first that \ref{dyncon} is true. Let us
assume by contradiction that there exists    an open entourage
$\bu\in\s{Ent}X$ for which $G_\bu/H$ is an infinite set.  Then
there exists a sequence $\{x_n, y_n\}\subset g_nS$ such that
$<x_n, y_n> \in\bv,$ where $\bv$ is the closed complement of $\bu$
in $\Th^2 X.$ Up to passing to a subsequence we obtain $x_n\to x$
and $y_n\to y$ and $<x,y>\in \bv.$ So we can choose closed
disjoint neighborhoods $L$ and $K$ of the points $x$ and $y$  such
that $g_nS\cap L\ne\emptyset$ and $g_nS\cap K\ne\emptyset$ for
infinitely many $n$ what is impossible.

  If, conversely  $K$ and $L$ are
disjoint closed subsets of $X$, let $\bu{=} (K\times L)'$. Then
the Definition implies that the set \ref{dyncon}  is at most
finite.

\medskip

b) Note that the definitions remain equivalent if one restricts to
the entourages of the form $\bold
u_\varepsilon{\leftrightharpoons}\{\{\g{p,q\}:\d(p,q}){<}\varepsilon\}$
($\varepsilon{>}0$) where $\d$ is a metric determining the
topology of $X$. In our case $\bold\Theta^2X/G$ is compact hence
we can restrict ourselves to the entourages that belong to a fixed
$G$-orbit that generates the filter $\s{Ent}X$ [Ge09, Prop E].
\setcounter{prop}2%
\begin{prop}\label{vistodyn} Let $G$ act 3-discontinuously on a
compactum $\widetilde T.$ An orbit $F{\subset}A$ of a subgroup $H$
of $G$ is visibly quasiconvex if and only if $H$ is dynamically
quasiconvex in $G$.
\end{prop}
\it Proof\rm. Since $H$ acts on $\widetilde T$ 3-discontinuously
we have $\partial F{=}\bold\Lambda H$. For $g{\in}G$ the set $gF$
is an orbit of the group $gHg^{-1}$. Thus
$\partial(gF){=}\bold\Lambda(gHg^{-1}){=}g\bold\Lambda H$. So, the
dynamical quasiconvexity of $H$ is equivalent to the \bf right
\rm$H$-finiteness of the sets of the form
$\{g{\in}G:\partial(gF){\notin}\s{Small}(\bold u)\}$ ($\bold
u{\in}\s{Ent}\widetilde T$).

We fix a reference vertex $a{\in}A$ and consider the generating
set of entourages of the form $\bold u_\varepsilon$ with respect
to the metric $\overline\d_a$. Since
$\overline\d_a(gx,gy){=}\overline\d_{g^{-1}a}(x,y)$, the following
property is equivalent to the dynamical quasiconvexity of $H$:

\hfil\penalty-10000 ($*$) for every $\varepsilon{>}0$ the set
$\mathcal
G_\varepsilon{\leftrightharpoons}\{g{\in}G:\s{diam}_{\overline\d_{ga}}
\partial F{\geqslant}\varepsilon\}$ is \bf left \rm$H$-finite.

 \medskip

 Since the left $H$-action preserves the $\s d$-distance from
$F$, if $\mathcal G_\varepsilon/H$ is finite then $\s d(\mathcal
G_\varepsilon a,F)$ is bounded. On the other hand, if $\s
d(\mathcal G_\varepsilon a,F)$ is bounded then $\mathcal
G_\varepsilon a$ is $H$-finite. Since the action
$G{\curvearrowright}A$ is properly discontinuous the set $\mathcal
G_\varepsilon$ is also $H$-finite. So ($*$) is equivalent to:

\hfil\penalty-10000($**$) $\s d(\mathcal G_\varepsilon a,F)$ is
bounded for every $\varepsilon{>}0$.

\medskip

Thus if $F$ is visibly quasiconvex then (**) is true for every
$a{\in} A$, so $H$ is dynamically quasiconvex.

Conversely, suppose that $\s d(\mathcal G_\delta
a,F){\leqslant}R_\delta$ for every $\delta{>}0$. Let $S\subset A$
be a finite set containing $a$ and intersecting each $G$-orbit in
$A$. Since change of the reference point is a bilipschitz
transformation, the ratio $\overline\d_x/\overline\d_y$
($x,y{\in}S$) is bounded.

Let $v{\in}\s V_\varepsilon F$, i.e,
$\s{diam}_{\overline\d_v}F{\geqslant}\varepsilon$. Then $v{\in}gS$
for some $g{\in}G$. So
$\s{diam}_{\overline\d_{ga}}F\geqslant\delta{\leftrightharpoons}
{\varepsilon\over C}$ for some uniform constant  $C$. We have
$g{\in}\mathcal G_\delta$, $\s d(ga,F){\leqslant}R_\delta$, $\s
d(v,F)\leqslant R_\delta{+}\s{diam_d}S$. So $F$ is visibly
quasiconvex.\qed

\subsection{Horocycles}

\bf Definition\rm. A bi-infinite $\alpha$-distorted path $\gamma :
\Bbb Z \to A$ is called $\alpha$\it-horocycle \rm  at $\g {p}\in
T$ if $\displaystyle\lim_{n\to\pm\infty}\gamma(n){=}\g p$.  We
call the unique limit point $\g p$ of $\gamma$      \it base \rm
of the horocycle.

Recall that  a limit point $x\in \bold\Lambda G$ is called {\it
conical} if there exists an infinite sequence of distinct elements
$g_n\in G$ and distinct points $a, b\in S $ such that $g_n(y)\to
a$ for all $y{\not=} x$ and  $g_n(x)\to b.$
\begin{prop}\label{horoc}
There is no $\alpha$-horocycle at conical point.
\end{prop}
 The proof  of this fact for quasigeodesic horocycles \cite[Lemma
3.6]{GP09}   works for $\alpha$-horocycles too.\qed

%\note{A slightly different proof, using the notion of `cone' from \cite[subsection 3.5]{Ge10}.\hfil\penalty-10000
%\hfil\penalty-10000
%Let $S$ be a cone in $G$ with vertex $\g p$.
%We have $\overline{S\g p}{\cap}\partial_1S{=}\varnothing$.
%By \ref{hullClosed}, the set $F{\leftrightharpoons}\s H_\alpha(\overline{S\g p})$ is closed in $\widetilde T$.
%It is clearly disjoint from $\partial_1S$.
%For an $\alpha$-horocycle $\gamma$ at $\g p$ and $v{\in}A{\cap}\s{Im}\gamma$
%the set $W{\leftrightharpoons}v{\times}F$ is closed in $\widetilde T^2$ and does not intersect
%the limit crosses for $S$. Hence there exists $s{\in}S$ whose graph does not meet $W$.
%On the other hand, $(v,s(v)){\in}W$. A contradiction.\qed}
\begin{prop}\label{2horo} Suppose that the action $G{\curvearrowright}\widetilde T$
is  3-discontinuous and 2-cocompact. Then there exists
$\varepsilon{>}0$ such that if $\alpha$-horocycles $\gamma,\delta$
with {\bf distinct} bases $\mathfrak{p,q}$ meet $a{\in}A$ then
$\overline\d_a(\g{p,q})\geqslant\varepsilon$.
\end{prop}
\bf Remark\rm. This statement could be easily deduced from the
results of \cite{Ge09}. However this should require the theory of
linkness and   betweenness relation developed in \cite{Ge09}. We
prefer to give a simple independent proof. Here for the  first
time we  use the 2-cocompactness of the action
$G{\curvearrowright}\widetilde T$.

\medskip

\it Proof\rm. By 2-cocompactness there exists $\varrho{>}0$ such
that for every different $x,y{\in}T$ one has
$\overline\d_v(x,y){>}\varrho$ for some $v{\in}A$. Let $v$ be such
a vertex for $\g p$ and $\g q$. Then the vertex $a$ does not
belong to at least one of the sets $\s
N^{\overline\d_v}_{\varrho/2}\g p$, $\s
N^{\overline\d_v}_{\varrho/2}\g q$. So by \ref{hullLoc} $\s
d(v,a){\leqslant}s{\leftrightharpoons}s(\varrho/2,\alpha)$
 and $\overline\d_a(\g{p,q}){>}\lambda^s\varrho$.\qed
\vskip5pt \bf Corollary\sl.\ Every $a\in A$ can belong to a
uniformly bounded number of $\alpha$-horocycles at different
bases.

\it Proof\rm. Since $A$ is $G$-finite it is enough to prove that
every $a\in A$ can belong to at most finitely many
$\alpha$-horocycles with different bases. Suppose not and $a\in
\bigcap_{i\in I}\gamma_i$ where $\gamma_i$ is an
$\alpha$-horocycle at $p_i$ and  $\vert I\vert{=}\infty.$  Since
$\widetilde T$ is a compactum the infinite set $P{=}\{p_i\ \vert\
i\in I\}$   must contain a convergent subsequence which is
impossible by Proposition \ref{2horo} .\qed
\section{Horospheres}
\subsection{Systems of horospheres}\label{systHor}
Let $\s{St}_Ga$ denote  the stabilizer $\{g{\in} G : ga{=}a\}$ of
a point $a\in \widetilde T$ in $G.$ We make use of the following
obvious property of the actions of a group $G$ on sets:

\begin{prop}\label{gFinProd} For $G$-finite $G$-sets $A,B$ the
following properties of a $G$-set $S{\subset} A{\times}B$ are
equivalent:\hfil\penalty-10000 $\s a:$ $S$ is
$G$-finite;\hfil\penalty-10000 $\s b:$ for every $a{\in}A$ the set
$S\cap\{a\}{\times}B$ is $\s{St}_Ga$-finite;\hfil\penalty-10000
$\s c:$ for every $b{\in}B$ the set $S\cap A{\times}\{b\}$ is
$\s{St}_Gb$-finite.\end{prop}\qed

We apply \ref{gFinProd} to the case when $A$ is as above and $B{\leftrightharpoons}\s{Par}$ is the set of parabolic points.
Taking into account that $\s{St}_Ga$ is finite for each $a{\in}A$ we have the following corollary:

\begin{prop}\label{horosph} The following properties of $G$-set $S{\subset} A{\times}\s{Par}$
are equivalent:\hfil\penalty-10000 $\s a:$ $S$ is
$G$-finite;\hfil\penalty-10000 $\s b:$ for every $a{\in}A$ the set
$S\cap\{a\}{\times}\s{Par}$ is finite;\hfil\penalty-10000 $\s c:$
for every $\g p{\in}\s{Par}$ the set $S\cap A{\times}\{\g p\}$ is
$\s{St}_G\g p$-finite.\end{prop}\qed

\bf Definition\rm. Any $G$-invariant $G$-finite subset $S$ of
$A{\times}\s{Par}$ determines a \it system of horospheres\rm. For
such   $S$ each set $S_{\g p}{\leftrightharpoons}\{a{\in}A:(a,\g
p){\in}S\}$ is called horosphere \it at \rm the parabolic point
$\g p$. The entire set $S$ is completely defined by the family
$\s{Par}{\ni}\g p\mapsto S_{\g p}{\subset}A$. So such a family
also \it determines a system of horospheres\rm. To satisfy the
conditions of \ref{horosph} this map should be $G$-equivariant and
each $S_{\g p}$ should be $\s{St}_G\g p$-finite.

\medskip

Examples:\hfil\penalty-10000 1. The set $\{(\mathbf a,\g
p):\mathbf a\underset{A,k}\#\g p\}$ studied in
\cite[6.10--7.2]{Ge09} for fixed $k{\geqslant}2$ determines a
system of horospheres. \hfil\penalty-10000 2. (most important for
this paper) For a distortion function $\alpha$, the family $\g
p\mapsto\s H_\alpha\g p$ is a system of horospheres. The condition
$(\s b)$ of \ref{horosph} follows from Corollary of \ref{2horo}.
This family has been studied in \cite{GP09} for affine functions
$\alpha$.\hfil\penalty-10000 3. If $\g p\mapsto S_{\g p}$ is a
system of horospheres and $r$ is a positive integer then the
family $\g p\mapsto\s N_r^{\s d}S_{\g p}$ is also a system of
horospheres.\hfil\penalty-10000 4. If $\g p\mapsto S_{\g p}$ is a
system of horospheres and $\alpha$ is an appropriate distortion
function then the family $\g p\mapsto\s H_\alpha S_{\g p}$ is also
a system of horospheres.\hfil\penalty-10000 5. The union of two
systems of horospheres is a system of horospheres.
\begin{prop}\label{nearHoro}
Let $S$ be a system of horospheres. Then for every
$r{{\geqslant}0}$ the set\hfil\penalty-10000
$\{(\g{p,q}){\in}\s{Par}^2:\s d(S_{\g p},S_{\g
q}){\leqslant}r\}\}$ is $G$-finite.
\end{prop}
\it Proof\rm. By passing to the system of horospheres $\g
p\mapsto\s N_rS_{\g p}$ the problem reduces to the case $r{=}0$.
In this case consider the set
$\{(a,\g{p,q}){\in}A{\times}\s{Par}^2:a{\in}S_{\g p}{\cap}S_{\g
q}\}$. It  is $G$-finite   by \ref{horosph}. The set
$\{(\g{p,q}){\in}\s{Par}^2:S_{\g p}{\cap}S_{\g
q}{\ne}\varnothing\}$ is the image of the latter one by the
$G$-equivariant map of forgetting the $a$-component. So it is also
$G$-finite.\qed
\begin{prop}\label{horHor}
Given a system of horospheres $\g p\mapsto S_{\g p}$ there exists
a positive $C$ such that\hfil\penalty-10000 $\s{diam_dPr}_{S_{\g
p}}S_{\g q}\leqslant C$ for each pair $\{\g{p,q}\}$ of distinct
parabolic points.
\end{prop}
\it Proof\rm. Let $\g q{\in}\s{Par}$. Since the action $\s{St}_G\g
q{\curvearrowright}(\widetilde T\setminus \g q)$ is cocompact, the
subgroup $S_{\g q}{=}\s{St}_G\g q$ is quasiconvex by
\ref{subgrConvex} and hence visibly quasiconvex by Corollary of
\ref{eHull}. Since $\s{Par}/G$ is finite, the visible
quasiconvexity function $\s Q_{S_{\g q}}$ (see \ref{eConvex}) can
be chosen independently of $\g q$. We denote any such function by
$\s Q_S$. That is: $\forall\g
q{\in}\s{Par}\forall\varepsilon{>}0\exists r{\leftrightharpoons}
\s Q_S(\varepsilon):\{a{\in}A:\s{diam}_{\overline\d_a}S_{\g
q}{\geqslant}\varepsilon\}\subset\s N_rS_{\g q}$.

Since $\s{Par}/G$ is finite it suffices to find $C$ for a
particular $\g p$. We thus fix it and denote
$H{\leftrightharpoons}\s{St}_G\g p$,
$\Sigma{\leftrightharpoons}S_{\g p}$.

Let $K$ be a compact fundamental set for
$H{\curvearrowright}(T\setminus \g p)$. So $K{\cap}{\overline
\Sigma}{=}\varnothing$. Since ${\overline \Sigma}$ is weakly
homogeneous by \ref{prOfA} the set $\s{Pr}_{\Sigma}K$ is finite
and
$\varrho{\leftrightharpoons}\s{min}\{\overline\d_v(\Sigma,K):v{\in}\s{Pr}_{\Sigma}K\}>0$.
By \ref{nearHoro} the set $P{\leftrightharpoons} \{\g
q{\in}\s{Par}:\s d(\Sigma,S_{\g
q}){\leqslant}r{\leftrightharpoons}\s Q_S(\varrho/2)\}$ is
$H$-finite. So
$C_1{\leftrightharpoons}\s{sup\{diam_d(Pr}_{\Sigma}S_{\g q}):\g
q{\in}P\}<\infty$.

If now $\g q{\not\in}K\setminus P$ then  up to applying an element
from $H$ we can assume that ${\g q}{\in} K$.   For $v\in
\s{Pr}_{\Sigma}\g q$
 we have $d(v, S_{\g q}) > r$ and $\s{diam}_{\overline \d_v}(S_{\g q}){<}\varrho/2.$ Thus
 $$\overline\d_v(\g p,S_{\g q}){\geq}\overline\d_v(\g p, \g q) -
 \overline\d_v(\g q, S_{\g q}){\geq}\overline\d_v(\Sigma, K)-
 \s{diam}_{\overline \d_v}(S_{\g q}){\geq}
\varrho{-}(\varrho/2){=}\varrho/2.$$ Hence for the number
$\rho(S_{\g q},\Sigma)$, defined in  \ref{numbrho}, we have
$\rho(S_{\g q},\Sigma){\geq}\varrho/2$. By \ref{prUniSize} we
obtain $C_2{\leftrightharpoons}\s{sup\{diam_d(Pr}_{\Sigma}S_{\g
q}):\g q{\notin}P\}<\infty$.  So we put
$C{\leftrightharpoons}\s{max}\{C_1,C_2\}$.\qed \vskip5pt \bf
Corollary\sl. Given a system $S$ of horospheres there exists a
positive number $C$ such that

$\s{diam_d}(S_{\g p}{\cap}S_{\g q}){\leqslant}C$ for each pair
$\{\g{p,q}\}$ of distinct parabolic points\rm.\qed
\subsection{Horospherical depth}\label{depth}
\def\h{horospherical }%
\bf Definition\rm. Let $\alpha$ be a distortion function and let
$\gamma:I\to A$ be a path. For $i{\in}I$ we define the \it\h depth
\rm of $i$ as \setcounter{equation}0
\begin{equation}\s{depth}_\alpha(i,\gamma)\leftrightharpoons\s{sup}\{r{\in}\Bbb N:\s N_ri{\subset}I\text{ and }
\exists\g p{\in}\s{Par}\ \gamma(\s N_ri)\subset\s H_\alpha\g
p\}.\end{equation} To take into account multiple points
 we put   for $v{\in}\s{Im}\gamma$
\begin{equation}\s{depth}_\alpha(v,\gamma)
{\leftrightharpoons}\s{inf}\{\s{depth}_\alpha(i,\gamma):\gamma(i){=}v\}.\end{equation}

\medskip

Applying the above Corollary to the system of the horospheres $\g
p\mapsto\s H_\alpha\g p,$ we obtain that there exists a constant
$h$ such that if $\gamma$ is $\alpha$-distorted and
$\s{depth}_\alpha(i,\gamma){\geqslant}h$ then there is exactly one
$\g p{\in}\s{Par}$ such that $\gamma(\s N_hi){\subset}\s
H_\alpha\g p$. We call such $h$ the \it critical depth value \rm
for $\alpha$.

\medskip

%\bf Definition\rm. Let $\gamma:I\to A$ be an $\alpha$-distorted
%path. A point $i{\in}I$ is said to be $(\alpha)$\it-\h\rm for
%$\gamma$ if for the critical depth value $h$ the $h$-neighborhood
%$\s N_ei{\leftrightharpoons}[i{-}h,i{+}h]$ is contained in $I$ and
%$\gamma(\s N_hi){\subset}\s H_\alpha\g p$ for some $\g
%p{\in}\s{Par}$. Note that this $\g p$ is unique if exists by the
%choice of $h$. A point that is not \h is called ($\alpha$-)\it
%non-horospherical\rm. When $\alpha$ is fixed we omit the prefix
%`$\alpha$-'.

%\medskip

%We also adopt the following convention: if a path $\gamma$ is
%injective and $i$ is a (non-)\h point for $\gamma$ we say that
%$\gamma(i)$ is a (non-)\h vertex of $\gamma$.

\vskip3pt Until the end of this subsection we fix an appropriate
distortion function $\alpha$.

For a vertex $v{\in}A$ denote by $\s{NH}_{v,e,\alpha}$ the set of
\bf finite \rm $\alpha$-distorted paths $\gamma$ of length
${>}\alpha_0{=}\alpha(1)$ such that
$\s{depth}_\alpha(v,\gamma){\leqslant}e$. Note that
$\partial\gamma$ is a proper pair for
$\gamma{\in}\s{NH}_{v,e,\alpha}$. \setcounter{prop}2
\begin{prop}\label{aEnt}
The set $\{\partial\gamma:\gamma{\in}\s{NH}_{v,e,\alpha}\}$ is bounded in $\bold\Theta^2\widetilde T$.
\end{prop}
\it Proof\rm. Otherwise there is a limit point $\{\s{p,p}\}$ for
this set. Since $A$ is discrete we have $\g p{\in}T$. By
compactness of the Tikhonoff topology there exists an
$\alpha$-horocycle $\gamma:\Bbb Z\to A$ at $\g p$ such that
$\gamma(0){=}v$ and for every finite segment $I{\subset}\Bbb Z$
there exists $\delta{\in}\s{NH}_{v,e,\alpha}$ such that
$\gamma|_I{=}\delta|_I$.

We have $\s{Im}\gamma{\subset}\s H_\alpha\g p$ hence
$\s{depth}_\alpha(0,\gamma){=}\infty$ contradicting with the
boundness of $\s{depth}_\alpha(0,*)$ on $\s{NH}_{v,e,\alpha}$.\qed

 \vskip5pt \bf
Corollary\sl. $\exists\varepsilon{>}0 \ \forall v{\in} A \ \forall
\gamma{\in}\s{NH}_{v,e,\alpha}\ :\
\s{diam}_{\overline\d_v}(\partial\gamma){>}\varepsilon$.\rm\qed

\section{Relative geodesics}
\subsection{Lifts}\label{lifts}
Let $S$ be a system of horospheres. We attach to our graph
$\Gamma$ new edges joining by an edge of length 1 each pair of
points that belong to an horosphere. The new graph is called \it
relative graph \rm and is denoted by $\Delta$. The corresponding
\it relative distance \rm function is denoted by $\overline{\s
d}$. The edges of $\Delta^1\setminus \Gamma^1$ are called \it
\h\rm and those belonging to $S_{\g p}$ are called $\g
p$\it-horospherical\rm. A change of the system of horospheres
yields a quasi-isometry of the relative graphs. To distinguish
pathes in $\Gamma$ and $\Delta$ we speak of $\Gamma$\it-paths \rm
and $\Delta$\it-paths\rm.

A $\Gamma$-path $\gamma$ is called a \it lift \rm of a
$\Delta$-path $\delta$ if these pathes have the same non-\h edges,
and, instead of any \h edge of $\delta$  in $\gamma$ one has a $\s
d$-geodesic segment with the same endpoints. Every subpath
$\gamma\vert_I$ (and the interval $I$) of $\gamma$ coming from an
edge of $\delta$ we call $\delta$\it-piece \rm of $\gamma$.  So,
to each edge of $\delta$ (called $\delta$-edge), \h or not, there
corresponds exactly one $\delta$-piece of $\gamma$. Note that a
lift of a $\overline{\s d}$-geodesic $\Delta$-path is not
necessarily injective.

A subpath $\gamma\vert_I$ (and the interval $I$) of $\gamma$ is
said to be \it integral \rm if it is a lift of some subpath of
$\delta$.
\begin{prop}\label{liftUndistorted}
There exists a quadratic polynomial $\alpha$ such that any lift $\gamma$ of any $\overline{\s d}$-geodesic $\Delta$-path $\delta$
is $\alpha$-distorted.
 Moreover, $\s{depth}_\alpha(v,\gamma)$ is uniformely bounded  for every
 $v{\in}\s{Im}\delta.$
\end{prop}
\it Proof\rm. Consider a lift $\gamma$ of a $\overline{\s
d}$-geodesic path $\delta:[j_-,j_+]\to A$. Let $I{=}[i_-,i_+]$ be
the corresponding  subinterval of $\s{Dom}\gamma$. We must prove
that $\s{diam}I{\leqslant}\alpha_n$ where
$n{\leftrightharpoons}\s{diam_d}\gamma(\partial I)$ and $\alpha$
is a quadratic polynomial that does not depend on $\gamma$ and
$\delta$.

Denote $P{\leftrightharpoons}\{\g p{\in}\s{Par}:$ there is a $\g
p$-\h edge in $\delta\}$. Since $\delta$ is $\overline{\s
d}$-geodesic it has exactly one $\g p$-\h edge for each $\g
p{\in}P$. Denote by $\gamma_{\g p}:I_{\g p}\to A$ the
corresponding $\s d$-geodesic segment of $\gamma$. Note that
$\s{Im}\gamma_{\g p}\subset\s{H_{id}}S_{\g p}$. By the
quasiconvexity of horospheres (see \ref{subgrConvex}) there exists
$r$ such that $\s{H_{id}}S_{\g p}\subset\s N_rS_{\g p}$ for each
$\g p{\in}\s{Par}$. We fix such $r$. Let $C_1$ be a maximum of the
constants determined in Proposition \ref{qIso} for
$E\leftrightharpoons\s N_rS_{\g p}$ ($\g p{\in}\s{Par}$). Let
$C_2$ be the constant determined by \ref{horHor} for the system of
horospheres $\g p\mapsto\s N_rS_{\g p}$.

We can assume that the points $i_-$ and $i_+$ belong respectively
to the first and to the last intervals of the set $\{I_{\g p}:\g
p{\in}P\}$. So $\overline{\s d} (\delta(j_\pm),\gamma(i_\pm))
\leqslant r{+}1$. Since $\s{length_{\overline d}\delta {=}
diam_{\overline
d}}\partial\delta\leqslant2r{+}2{+}\s{diam_d}\partial\gamma{=}n{+}2r{+}2$
we have $0\leqslant r_1\leftrightharpoons|\{$non-\h edges of
$\delta\}|\leqslant n{+}2r{+}2{-}|P|$.

 We now claim that $\s{diam}(I{\cap}I_{\g p})\leqslant C(2n{+}2r{+}1)$ where $C{\leftrightharpoons}\s{max}\{C_1,C_2\}$.
Indeed let $\beta$ be a $\s d$-geodesic segment between the
endpoints of $\gamma$. Consider the path $\omega$ that joins the
endpoints of $\gamma(\partial(I{\cap}I_{\g p}))$ formed by $\beta$
and the two pieces of $\gamma$ between the endpoints of $\gamma$
and the corresponding endpoints of $\gamma_{\g p}$ (one of the
pieces can be empty). We then ``project'' $\omega$ onto $\s
N_rS_{\g p}$ as follows: for each vertex
$v\in\s{Im\beta{\cup}Im}\delta{\cup}\partial\gamma$ we choose a
vertex in $\s{Pr}_{\s N_rS_{\g p}}v$ and join them  by $\s
d$-geodesic segments. Each edge of $\beta$ and each
non-horospherical edge of $\delta$ gives at most $C_1$ edges of
the projection. Each piece $\gamma_{\g q}$ ($\g q{\in}P\setminus
\g p$), corresponding to an \h edge of $\delta$
  gives at most $C_2$ edges in the projection. The curve
   $\omega$ does not contain
  $\g p$-horospherical edges. Thus the $\s
d$-distance between the endpoints of $\gamma|_{I{\cap}I_{\g p}}$
is at most $C_1{\cdot}(n{+}r_1)+C_2(|P|{-}1)\leqslant
C(2n{+}2r{+}1)$ and our claim is proved.

Since $\delta$ has at most $n{+}2r{+}2$ edges (either \h or not)
we have the following estimate\hfil\penalty-10000
$\s{length_d}\gamma\leqslant
C(n{+}2r{+}2)(2n{+}2r{+}1){=}\alpha_n$ where $\alpha$ is a
polynomial of degree 2.  Thus the lift of a $\Delta$-geodesic path
$\delta$ is an $\alpha$-distorted path $\gamma$ in $\Gamma$
proving the first part of the Proposition. \vskip3pt To estimate
the $\alpha$-\h depth
%\note{it seems that the rest of the proof is valid for each $\alpha$. Is it reasonable
%to separate it into another proposition?}
of the vertices of $\delta$ in $\gamma$ we fix a number $s$ (see
\ref{subgrConvex}) such that $\s H_\alpha\g p{\subset}\s H_\alpha
S_{\g p}\subset\s N_s^{\s d}S_{\g p}$ for every $\g
p{\in}\s{Par}$.

Let $v{\in}\s{Im}\delta$. Assume that $\gamma(0){=}v$.

Let $K{=}[k_-,k_+]$ be a maximal subinterval of
$I{\leftrightharpoons}\s{Dom}\gamma$ containing $0$ such that
$\gamma(\partial K){\subset}\s H_\alpha\g p$ for some $\g p$. We
fix such $\g p$.

We have $\s{depth}_\alpha(0,\gamma){=}\s{min}\{|k_-|,|k_+|\}$.
Since each two points of $\gamma K$ can be joined by a
$\overline{\s d}$-geodesic path of length at most $2s{+}1$ through
$S_{\g p}$ we have $\s{diam_{\overline d}}(\delta^{-1}\gamma
(K))\leqslant2s{+}1$.

Let $L{=}[l_-,l_+]$ be the largest integral subinterval of $K$ and
let $M{=}[m_-,m_+]$ be the smallest integral interval containing
$K$. So we have $0{\in}L{\subset}K{\subset}M{\subset}I$.

Since $\delta$ is $\overline{\s d}$-geodesic, at least one of the
integral intervals $M_-{\leftrightharpoons}[m_-,0]$,
$M_+{\leftrightharpoons}[0,m_+]$ does not contain a $S_{\g
p}$-edge of $\delta$. Let us assume that it is $M_+$. Note that
$k_+ > l_+$ only if $\displaystyle \gamma\vert_{[l_+,k_+]}$
belongs to $\s N_s S_{\g q}$ for some $\g q{\in}
\s{Par}\setminus\g p.$ Thus $k_+{-}l_+\leqslant
c{\leftrightharpoons} \s{max\{diam_d(N}_sS_{\g p}{\cap}\s N_sS_{\g
q}):\{\g{p,q}\}{\in}\bold\Theta^2\s{Par}\}$
 (see \ref{nearHoro}). Since each $\delta$-edge in $[0, l_+]$ yields at most $c$ edges in
 $\gamma\cap \s N_sS_{\g p}{\cap}\s N_sS_{\g q}$ we obtain
$l_+{\leqslant}(2s{+}1)c$. Hence
$\s{depth}_\alpha(0,\gamma)\leqslant(2s{+}2)c$.\qed

\subsection{Relative hull}\label{relHull}
For a set $F{\subset}A$ define its \it relative hull \rm$\s H_{\fam0rel}F{\leftrightharpoons}
{\cup}\{\s{Im}\delta:\delta$ is a $\s{\overline d}$-geodesic $\Delta$-path with $\partial\delta{\subset}F\}$.
A set is said to be \it relatively quasiconvex \rm if $\s H_{\fam0rel}F\subset\s N_rF$ for some $r{<}\infty$.

From now on we suppose that $\sum_{n\geqslant0}n^2f_n{<}\infty$
for our scaling function $f$ (\ref{scale}). For example we can
take $f_n{=}(n{+}1)^{-3-\varepsilon}$ for any $\varepsilon{>}0$.
Thus any pair $(f,\alpha)$ where $\alpha$ is a quadratic
polynomial is appropriate.

Denote by $\lambda$ the decay rate (\ref{scale}) of $f$.
\begin{prop}\label{relDyn}
There exists a function $r{=}r(\varepsilon)$ such that for every
$F{\subset}A$, $\s V_\varepsilon F{\subset}\s N_r^{\s d}\s
H_{\fam0rel}F$. In particular, every relatively quasiconvex set is
visibly quasiconvex (see \ref{eConvex}).
\end{prop}
\it Proof\rm. Let $\alpha$ be the distortion function from
\ref{liftUndistorted}. For $\varepsilon{>}0$ we will find $r$ that
depends only on $\alpha, \varepsilon$, $\lambda$, the Karlsson
functions, and the convergence functions (see \ref{prExtent}) of
the horospheres.

Let $a{\in}\s V_\varepsilon F$, i.e,
$\s{diam}_{\overline\d_a}F{\geqslant}\varepsilon$. Connect a pair
of points  of $\overline\d_a$-diameter
${\geqslant}{\varepsilon\over2}$ in $F$ with a $\overline{\s
d}$-geodesic $\Delta$-path $\delta$ and consider its lift
$\gamma:I{=}[i_-,i_+]\to A$. By Proposition \ref{liftUndistorted}
$\gamma$ is $\alpha$-distorted. Since the pair $(f,\alpha)$ is
appropriate we can assume that $\s d(a,b)\leqslant
d_1{\leftrightharpoons}\s K_\alpha({\varepsilon\over2})$
 where $b{\leftrightharpoons}\gamma(0)$.

 If $b\in \delta^0$ we are
 done, so suppose not.
We have
$\s{diam}_{\overline\d_b}\partial\gamma\geqslant\rho{\leftrightharpoons}\lambda^{d_1}{\varepsilon\over2}$.
Let $J{=}[j_-,j_+]$ be the $\delta$-piece (see \ref{lifts}) of $I$
containing $0$ and let $J_-{=}[i_-,j_-]$, $J_+{=}[j_+,i_+]$ be the
complementary subintervals of $I$. By the $\triangle$-inequality,
at least one of the numbers
$\overline\d_b(\gamma(j_+),\gamma(i_+))$,
$\overline\d_b(\gamma(j_-),\gamma(j_+))$,
$\overline\d_b(\gamma(i_-),\gamma(j_-))$ should be
${\geqslant}{\rho\over3}$. Respectively, consider these three
cases. The third one reduces obviously to the first.

In the first case we have $\s d(b,c)\leqslant
d_2{\leftrightharpoons}\s K_\alpha({\rho\over3})$ for some
$c{=}\gamma(k)$, $k{\in}J_+$. As  $\gamma$ is $\alpha$-distorted
  we have $k{\leqslant}\alpha(d_2)$. Since
$j_+{\in}[0,k]\cap\delta^0$ it follows that $\s
d(b,\s{Im}\delta){\leqslant}\alpha(d_2)$.

\hspace*{5cm}{\epsfxsize=3.5in\epsfbox{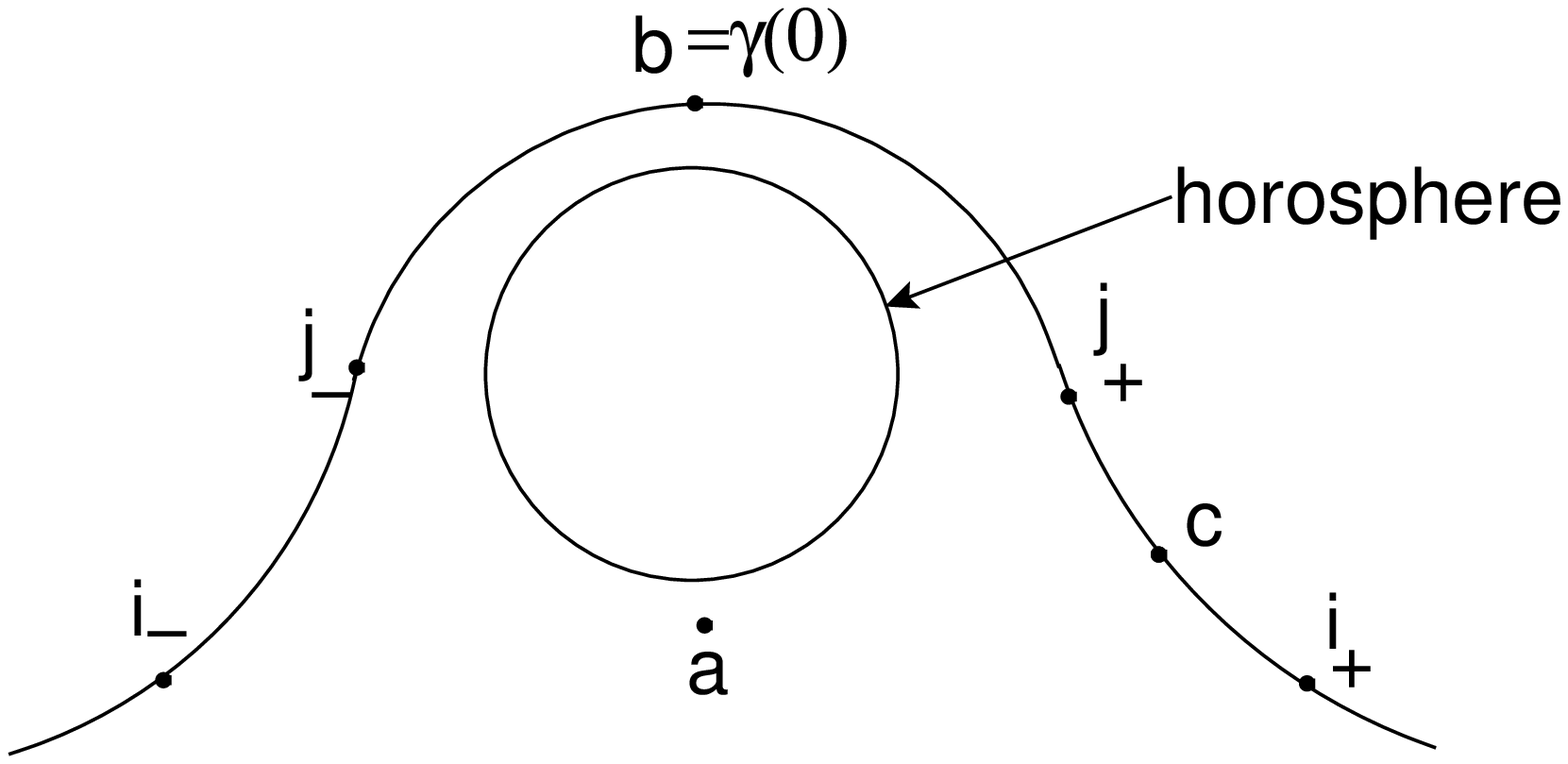}}

\hspace*{6.7cm}{Lift of a relative geodesic}

\bigskip

\bigskip

In the second case since $b\not\in\delta^0$ there exists an
horosphere $S_{\g p}\ (\g p{\in}\s{Par})$ such that
$\gamma(\partial J){\subset}S_{\g p}$, $\gamma|_J$ is geodesic,
and $b{\in}E{\leftrightharpoons}\s{H_{id}}S_{\g p}$. We claim that
$\s d(b,\s{Im}\delta)\leqslant d_3{\leftrightharpoons}\s
C_E({\rho\over6})$ (see \ref{CE}). If not, since
$\gamma(j_{\pm})\in \s{Im}\delta$, we have $d(b, \gamma(j_\pm))
{>}d_3$. By \ref{CEa} it follows that
$\overline\d_b(\gamma(j_\pm), p{=}\partial E){\leq}{\rho\over6}$
and so $\overline\d_b(\gamma(j_+), \gamma(j_-))\leq {\rho\over3}$
which is impossible. Since the set $S$ of the horospheres is
$G$-finite the above constant $d_3$ is uniform for every ${\g
p}\in \s{Par}.¤$

So, in either case we have the uniform bound

 $$\s d(a,\s
H_{\fam0rel}F)\leq \s d(a,b)+\s d(b, \s H_{\fam0rel}F) \leq d_1 +
\max\{d_3,\alpha(d_2)\}$$ as claimed.\qed

\medskip

\bf Definition\rm. Let $\alpha$ be a distortion function and let
$e$ be a positive integer. Define the $(\alpha,e)$\it-hull \rm of
a set $F{\subset}A$ as

\rm$\s H_{\alpha,e}F{\leftrightharpoons}\{\gamma(i):\gamma $ is an
$\alpha$-distorted path with $\partial\gamma{\subset}F$ and
$\s{depth}_\alpha(i,\gamma){\leqslant}e\}$.

A set $F{\subset}A$ is said to be \it relatively
$\alpha$-quasiconvex \rm if   its $(\alpha,e)$-hull is within a
bounded distance from $F$.

\medskip

It follows from the corollary of \ref{horHor} that this notion of
quasiconvexity does not depend on $e$ when $e$ is sufficiently
large.
\begin{prop}\label{deltaAlpha} There exists a function $\varepsilon{=}\varepsilon(\alpha,e)$
such that $\s H_{\alpha,e}F{\subset}\s V_\varepsilon F$ for every
$F{\subset}A$. In particular, every visibly quasiconvex set is
$\alpha$-quasiconvex for any appropriate distortion function
$\alpha$.\end{prop} \it Proof\rm. For $v{\in}A$ the number
$\s{inf\{diam}_{\overline\d_v}\partial\gamma:\gamma{\in}\s{NH}_{v,e,\alpha}\}$
is positive by the Corollary of \ref{aEnt}. Since $A$ is
$G$-finite the number $\varepsilon{\leftrightharpoons}
\s{inf\{diam}_{\overline\d_v}\partial\gamma:v{\in}A,\gamma{\in}\s{NH}_{v,e,\alpha}\}$
is also positive. It follows that if $v{\in}\s H_{\alpha,e}F$ then
$v{\in}\s V_\varepsilon F$.\qed
\begin{prop}Let $\alpha$ be a distortion function from \ref{liftUndistorted}.
Then there exists a number $v$ such that $\s
H_{\fam0rel}F{\subset}\s H_{\alpha,e}F$ for arbitrary set
$F{\subset}A$. In particular every relative $\alpha$-quasiconvex
set $F$ is relatively quasiconvex.\end{prop} \it Proof\rm. The
result follows immediately from \ref{liftUndistorted}.\qed
\vskip3pt Summing up the results of Subsection \ref{relHull} and
Proposition \ref{vistodyn}, we obtain. \vskip5pt \bf Theorem A\sl.
Let a finitely generated discrete group $G$ act 3-discontinuously
and 2-cocompactly on a compactum $\widetilde T$. The following
properties of a subset $F$ of the discontinuity domain of the
action are equivalent:\hfil\penalty-10000
--- $F$ is relatively quasiconvex;\hfil\penalty-10000
--- $F$ is visibly quasiconvex;\hfil\penalty-10000
--- $F$ is relatively $\alpha$-quasiconvex where $\alpha$ is a quadratic polynomial with big enough coefficients.

Moreover, if $H$ is a subgroup of $G$ and $F$ is $H$-finite then the visible quasiconvexity of $F$ is
equivalent to the dynamical quasiconvexity of $H$ with respect to the action $G{\curvearrowright}\widetilde T$\rm.\qed
\section{The lifts of geodesics from the relative graph and some applications.}
\subsection{Hyperbolicity of the relative graph}
As one of the applications of our methods we give an easy proof of
the main results of \cite{Ya04}. We first suppose that $G$ is a
finitely generated relatively hyperbolic group admitting a
geometrically finite convergence action $G\curvearrowright T$, or
equivalently the action is $2$-cocompact. Let $\Gamma$ be a
locally finite, connected, $G$-graph. Denote by $\Delta$ the
corresponding relative graph with respect to the system of
horospheres (see 6.1).

Our first aim is to show that the relative graph is Gromov
hyperbolic.

\begin{prop}\label{relthin}
There exists a constant $r$ such that, for every $\overline{\s
d}$-geodesic triangle in $\Delta$, every its side is within the
$r$-neighborhood in $\Delta$ of the union of the other two sides.
\end{prop}
\it Proof\rm. Let $\alpha$ be the distortion function from
\ref{liftUndistorted} and let $e$ be the upper bound for
$\s{depth}_\alpha(v,\gamma)$ from \ref{liftUndistorted}. Let
$\varepsilon{\leftrightharpoons}\varepsilon(\alpha,e)$ be the
number from \ref{deltaAlpha}.

Consider a $\overline{\s d}$-geodesic triangle with edges
$\delta,\delta',\delta''$. Let $\gamma,\gamma',\gamma''$ be the
lifts. We can assume that  $\delta(0){=}\gamma(0){=}v$. By
Proposition \ref{liftUndistorted} we have
$\gamma{\in}\s{NH}_{v,e,\alpha}$. Hence by Corollary of \ref{aEnt}
$\s{diam}_{\overline\d_v}\partial\gamma{\geqslant}\varepsilon$.
Thus one of the numbers $\s{diam}_{\overline\d_v}\partial\gamma'$,
$\s{diam}_{\overline\d_v}\partial\gamma''$ is
${\geqslant}{\varepsilon\over2}$. It follows
that\hfil\penalty-10000 $\s
d(v,\s{Im\gamma'{\cup}Im\gamma''){\leqslant}K}_\alpha({\varepsilon\over2})$.
So putting $r{=}\s K_\alpha({\varepsilon\over2}){+}1$ we obtain
 $\overline{\s d}(v, \s{Im\delta'{\cup}Im\delta'')}{\leqslant}r$.\qed

\medskip

\bf Remark\rm. One of the equivalent definitions of (strong)
relative hyperbolicity of a group was proposed by B. Bowditch. It
claims that a group is relatively hyperbolic if and only if it
possesses a cofinite action on a Gromov hyperbolic graph $\Delta$
(`cofinite' means that $\Delta^1$ is $G$-finite) which is \it
fine\rm, that is for every $n$ and every edge $e$ of $\Delta$ the
set of simple loops of length $n$ that path through $e$ is finite.

In our case the action $G{\curvearrowright}\Delta$ is not
cofinite, but the metric space $\Delta^0$ can be isometrically and
equivariantly embedded into a $G$-cofinite hyperbolic graph by the
following well-known construction: let $\widetilde\Delta$ be the
graph whose set of vertices is $A{\cup}\s{Par}$ and the set of
edges is $\Gamma^1{\cup}S$ (where $S$ is a $G$-finite subset of
horospheres in $A{\times}\s{Par}$ see 5.1). We consider on
$\widetilde\Delta^0$ the path-metric in which the $\Gamma$-edges
have length 1 and the $S$-edges have length $1\over2$. The
inclusion $\Delta^0{\hookrightarrow}\widetilde\Delta^0$ is an \bf
isometry \rm with respect to the path-metrics. Thus we can denote
the distance in $\widetilde\Delta$ by $\overline{\s d}$. Since
$A{=}\Delta^0{\subset}\widetilde\Delta^0{\subset}\s
N_{1/2}^{\overline{\s d}}A$, by Proposition \ref{relthin} the
graph $\widetilde\Delta$ is hyperbolic. The action
$G{\curvearrowright}\widetilde\Delta$ is cofinite. To prove the
finess of the graph $\widetilde\Delta$ we need the following lemma
motivated by \cite[Lemma 7.1]{Bo97}: \vskip3pt \bf Lemma\sl. There
exists a quadratic polynomial $\alpha$ such that for every simple
loop $\delta$ in $\widetilde\Delta$ and every its lift $\gamma$
one has
$$\s{length_d\gamma\leqslant\alpha(length_{\overline
d}}\delta).$$\rm

\it Proof\rm. Let $\delta$ be a simple loop in the graph
$\widetilde\Delta$ of length $n$. It can path at most once through
a parabolic point $\g p{\in}\s{Par}$. So we can suppose that
$\delta$ is a simple loop in the graph $\Delta$ having at most one
$\g p$-horospherical edge for each $\g p{\in}\s{Par}$. The
argument now repeats  the proof of the first part of
\ref{liftUndistorted} with obvious simplification. For the sake of
completeness we include it.

Let  $P{\leftrightharpoons}\{\g p{\in}\s{Par}:$ there is a $\g
p$-\h edge in $\delta\}$. Consider a lift $\gamma$ of $\delta$
(see \ref{lifts}). We can regard $\gamma$ as a map from the vertex
set $\Xi^0$ of a simplicial circle $\Xi$ taking edges to
$\Gamma$-edges. For every $\delta$-piece $\gamma_{\g p}$ denote by
$\omega_{\g p}$ the ``complementing path'' i.e. the restriction of
$\gamma$ onto $\Xi^0\setminus \s{Dom}\gamma_{\g p}$. By our
assumption $\omega_{\g p}$ does not path through $\gamma_{\g p}$
anymore. Thus $\omega_{\g p}$ consists of $n{-}1$ $\delta$-pieces
each piece is either $\gamma_{\g q}$ where $\g q{\in}P\setminus \g
p$ or a $\Gamma$-edge.

By projecting $\omega_{\g p}$ onto $S_{\g p}$ and comparing the
length of the resulting curve with the geodesic segment
$\gamma_{\g p}$ we have  $\s{length_d}\gamma_{\g p}\leqslant
C(n{-}1)$ where $C$ is the constant from the proof of
\ref{liftUndistorted}. Thus $\s{length_d}\gamma\leqslant
Cn(n{-}1)$.\qed

\medskip

\bf Corollary\sl \cite{Ya04}. For a finitely generated relatively
hyperbolic group $G$ the graph $\widetilde\Delta$ is fine.

\vskip3pt

\it Proof\rm. The graph $\Gamma$ is locally finite. So by the
above Lemma   there are at most finitely many lifts of a simple
loop of length $n$ in $\widetilde\Delta$ passing through a given
edge. It follows that $\widetilde\Delta$ is fine.\qed

 \medskip

The result of Yaman remains valid for relatively hyperbolic groups
without assuming their finite generatedness and even their
countability (cf with \cite{Hr10}).

\begin{prop}\label{genyam}
Let $G$ be a group acting 2-cocompactly and 3-discontinuously on a
compactum $T$. Then there exists a hyperbolic, $G$-cofinite graph
$\TD$ whose vertex stabilizers are all finite except the vertices
corresponding to the parabolic points for the action $G\act T.$
Furthermore the graph $\TD$ is fine.
\end{prop}

\it Proof\rm. We will use few facts from \cite{GP10}. The group
$G$ satisfying the above assumptions acts discontinuously on a
$G$-finite graph of entourages $\cal G.$ Denote by $P_i\
(i{=}1,...,n)$ the system of non-conjugate maximal parabolic
subgroups of $G$ for the action $G\act T.$ It is shown in [GP10,
Theorem A] that there is a graph $\widetilde{\cal G}$ obtained by
refinement of $\cal G$ such
 that all its connected components  are $G$-equivalent; and if
  $\G_0$ is a component of $\TG$ then its stabilizer   $G_0$
 is a finitely generated relatively hyperbolic subgroup   of
$G$ with respect to the system $Q_i{=}P_i\cap G_0$. The connected
components of $\widetilde{\cal G}$   are adjacent along the set of
parabolic points ${\g p}\in\s{Par}$ (not belonging to
$\widetilde{\cal G}$).

Let now $\TD$ be the graph  obtained  by joining every vertex of
$\ti\G$ belonging to an horosphere $S_{\g p}\in S$ with the
parabolic point $\g p$ by an edge of length $1\over 2$. The graph
$\TD$ is $G$-cofinite. Denote by $\TD_0$ the subgraph of $\TD$
corresponding to the component $\G_0$ of $\ti{\cal G}.$  By the
above Corollary  the graph $\TD_0$ is hyperbolic and fine.

There is an induced  action of $G$ on a bipartite graph $\T$ whose
vertices are of two types $\mathcal H$ and $\cal C$ corresponding
respectively to  the horospheres
 of $\widetilde{\cal G}$ (of {\it horospherical} type)
  and  to the connected components of
$\widetilde{\cal G}$ (of {\it non-horospherical} type). Two
vertices $H\in\cal H$ and $C\in \cal C$ are connected by an edge
in $\T$ if the corresponding horosphere $H$ and the component $C$
intersect. One can also obtain $\T$   from the graph $\TD$ by
contracting every component $g(\TD_0)\ (g\in G)$ into a vertex of
$\cal C$-type and every parabolic vertex ${\g p}\in\s{Par}$ into a
vertex of $\cal H$-type. By [GP10, Lemma 3.36] the graph $\T$ is a
tree. So every loop in $\TD$ is contained in $g(\TD_0)$ for some
$g\in G$. It follows that  the graph $\TD$ is itself a cofinite,
hyperbolic and fine.\qed

\subsection{The lifts of $\overline{\s d}$-geodesics are $\s d$-quasigeodesics}
Let $S$ be a system of horospheres. It follows from the definition
of a system of horospheres that the value $\s
C_S(\varepsilon){\leftrightharpoons}\s{sup\{C}_{S_{\g
p}}(\varepsilon):\g p{\in}\s{Par}\}$ (see \ref{prExtent}) is
finite for every $\varepsilon{>}0$.
\begin{prop}\label{hullT}
Given a system $S$ of horospheres there exists a number $d$ such
that if $\gamma:I\to A$ is $\alpha$-distorted, with
$\partial\gamma{\subset}S_{\g p}$ then $\gamma(I\setminus \s
N_d\partial I){\subset}\s H_\alpha\g p$.
\end{prop}
\it Proof\rm. If such $d$ were not exist one could find a sequence
of $\alpha$-distorted paths $\gamma_n:[i_n^-,i_n^+]\to A$ with
$\partial\gamma_n{\subset}S_{\g p_n}$, $\gamma_n(0){\notin}\s
H_\alpha\g p_n\ (\g p_n\in\s{Par})$, and $|i_n^\pm|\to\infty$.
Since the set $\s{Par}$ is $G$-finite by applying $G$ and passing
to a subsequence we can suppose that $\g p_n{=}\g p$. As
$\s{St}_G\g p$ acts cocompactly on $\widetilde T\setminus \g p$ we
can also assume that $\gamma_n(0){=}v$ do not depend on $n$. So by
passing to a  subsequence once more we can find a sequence of
paths that converges in the Tikhonoff topology to an infinite
$\alpha$-distorted path $\gamma$. Since $\partial S_{\g p}{=}\g
p$, it is an $\alpha$-horocycle and $\gamma(0){\notin}\s
H_\alpha\g p$. A contradiction.\qed

\medskip

Let $\alpha$ be a distortion function and let $e$ be a positive
integer.

\vskip3pt

\bf Definition\rm. An $\alpha$-distorted path $\gamma:I\to A$ is
called $e$\it-piecewise geodesic \rm if every subpath consisting
of points of $\alpha$-depth ${\geqslant}e$ (see \ref{depth}) is
geodesic.
\medskip

It follows from \ref{liftUndistorted} that for any system of
horospheres $S$ there exists $e$ such that every lift $\gamma$ of
a $\overline{\s d}$-geodesic path $\delta$ is $e$-piecewise
geodesic $\alpha$-distorted path for a quadratic polynomial
$\alpha_n.$

\begin{prop}\label{qGeo}
There exists a function $c{=}c(\alpha,e)$ such that every
$e$-piecewise geodesic $\alpha$-distorted path is
$\beta$-distorted where $\beta(n){=}cn{+}c\ (n\in\N)$.\end{prop}
\it Proof\rm. Consider a $e$-piecewise geodesic $\alpha$-distorted
path $\gamma:I{\leftrightharpoons}[0,i_+]\to A$ and a geodesic
path $\omega:J{\leftrightharpoons}[0,j_+]\to A$ with
$\gamma(0){=}\omega(0){=}a$, $\gamma(i_+){=}\omega(j_+){=}b$. Let
$h$ be the critical depth value (see \ref{depth}) for geodesics.
Denote\hfil\penalty-10000
$N{\leftrightharpoons}\{j{\in}J:\s{depth_{id}}(j,\omega){\leqslant}h\}$.
%(the points of $N$ are \it non-horospherical\rm)
By Corollary of \ref{aEnt} there exists $\varepsilon{>}0$ such
that, for $j{\in}N$, one has
$\overline\d_{\omega(j)}(a,b){\geqslant}\varepsilon$ and hence $\s
d(\omega(j),\s{Im}\gamma)\leqslant r{\leftrightharpoons}\s
K_\alpha(\varepsilon)$. Denote $s{\leftrightharpoons}\s
K_\alpha({\varepsilon\over2})$,
$t{\leftrightharpoons}2r{+}\alpha(r{+}s)$. \vskip3pt \bf Lemma\sl.
If $x,y{\in}N$ and $y{-}x>t$ and $\gamma(x_1){\in}\s
N_r(\omega(x))$, $\gamma(y_1){\in}\s N_r(\omega(y))$ then $x_1 <
y_1$\rm. \vskip3pt \it Proof\rm. Suppose not and $x_1 {\geqslant}
y_1.$  We have $d(\omega(x), [\gamma(y_1), \omega(y)]) >
t{-}r{=}r{+}\alpha(r{+}s){\geqslant}s$ where $[\gamma(y_1),
\omega(y)]\subset A$ is a geodesic of length at most $r$ between
$\gamma(y_1)$ and $\omega(y).$

By \ref{KaFn}
$\overline\d_{\omega(x)}(\omega(y),\gamma(y_1))\leqslant{\varepsilon\over2}$.
By Corollary of \ref{aEnt} we also have
$\overline\d_{\omega(x)}(a,\omega(y))\geqslant\varepsilon,$ so
$\overline\d_{\omega(x)}(a,\gamma(y_1))\geqslant{\varepsilon\over2}$.
Applying again \ref{KaFn} we obtain $\s
d(\omega(x),\gamma([0,y_1]))\leqslant s$.

 Let now
$x_2{\in}[0,y_1]$ be such that $\s
d(\omega(x),\gamma(x_2)){\leqslant}s$. Thus $ \s
d(\gamma(x_1),\gamma(x_2))\leqslant r{+}s$ and
$x_1{-}x_2\leqslant\alpha(r{+}s)$. Since $y_1{\in}[x_2,x_1]$ we
have

$\s d(\omega(x), \omega(y))\leqslant \s d(\omega(x),
\gamma(x_1)){+}\s d(\gamma(x_1), \gamma(y_1)){+}\s d (\gamma(y_1),
\omega(y))\leqslant 2r{+}\alpha(r{+}s){=}t$. Since $\omega$ is a
geodesic we obtain $y-x{\leqslant}t$. A contradiction.\qed
\vskip3pt We continue the proof of \ref{qGeo}. Subdivide the
interval $J$ into segments $J_k{\leftrightharpoons}[j_k,j_{k+1}]$
using the following inductive rule. Put
$j_0{\leftrightharpoons}0$. After the choice of $j_k$ if
$j_k{=}j_+$ then we finish. If not, define
$j_{k+1}{\leftrightharpoons}j_+$ if $j_+{-}j_k\leqslant t$.
Otherwise\hfil\penalty-10000
$j_{k+1}{\leftrightharpoons}\s{min}\{j{\in}N:j{>}j_k{+}t\}$.

Let $m$ be the biggest $k$ for which $j_k$ is defined.

 By the above argument, there exist $i_k{\in}[0,i_+]$
such that $\s d(\omega(j_k),\gamma(i_k))\leqslant r$ for $k{\in}[0,m{+}1]$.
By Lemma, the indices $i_k$ form an increasing
sequence. So $I$ gets subdivided into the segments $I_k{\leftrightharpoons}[i_k,i_{k+1}]$.
It suffices to find a linear polynomial $\beta$ such that
$\s{diam}I_k\leqslant\beta(\s{diam}J_k)$ for all $k$.

If $\s{diam}J_k\leqslant t{+}1$ then
$\s{diam}\partial(\gamma|_{I_k})\leqslant t{+}1{+}2r$
hence\hfil\penalty-10000 $\s{diam}I_k\leqslant\alpha(t{+}1{+}2r)
\leqslant\alpha(t{+}1{+}2r){\cdot} \s{diam}J_k$.

If $\s{diam}J_k>t{+}1$ then $J_k$ contains a piece of
$\s{id}$-depth ${>}h$. Hence, for a uniquely determined $\g
p{\in}\s{Par}$, the endpoints of $\omega|_{J_k}$ belong to the
$t{+}1$-neighborhood of $\s{H_{id}}\g p$.

Let $d$ be the constant from \ref{hullT} for the distortion function $\alpha$ and for the system of horospheres
$\g p\mapsto\s N_{t+1+r}\s{H_{id}}\g p$.

If $\s{diam}I_k>2(d{+}e{+}r{+}t)$ then, by \ref{hullT}, the
interval $I_k$ contains a nonempty subinterval
$I_k^{\fam0geo}{\leftrightharpoons}I_k\setminus \s N_{d+e}\partial
I_k\subset H_{\alpha}p$ of depth ${\geqslant}e$. By the hypothesis
$\gamma|_{I_k^{\fam0geo}}$ is a geodesic subpath. Thus
$\s{length_d}\gamma|_{I_k^{\fam0geo}}{=}\s{diam}I_k^{\fam0geo}\leqslant\s{diam}J_k+2(r{+}d{+}e)$.
Hence $\s{diam}I_k\leqslant\s{diam}J_k+2r{+}4d{+}4e$.

If $\s{diam}I_k\leqslant2(d{+}e{+}r{+}t)$ then also $\s{diam}I_k\leqslant2(d{+}e{+}r{+}t){\cdot}\s{diam}J_k$.

We have $[a, b]\subset \bigcup_k I_k$, so $b-a\leqslant
cd(\gamma(a), \gamma(b))+c$ where $c{=}\max\{2(d+e+r+t), 2r
+4d+4e, \alpha(t+1+2r)\}$. \qed

\vskip3pt

As a direct consequence of the above Proposition we obtain.

\vskip3pt

 \bf Corollary \cite[Thm. 1.12(4)]{DS05}. \sl The lift of
every $\overline{\s d}$-geodesic is $\s d$-quasigeodesic\rm.\qed

\section{Criteria for the subgroup quasiconvexity in RHG.}\setcounter{prop}1
\setcounter{equation}1

\subsection{Statement of the result.} The aim of this Section is to prove Theorem B giving
criteria for a subgroup of a relatively hyperbolic group to be
quasiconvex.

Let $q$ be a positive integer and let $\alpha$ be a distortion
function (see Subsection 2.2).

\medskip

\bf Definition\rm. A subset $E$ of a metric space $M$ is called
\it weakly $\alpha$-quasiconvex \rm if there is a positive integer
$q$ such that for each $x,y{\in}E$ there exists an
$\alpha$-distorted path $\gamma$ such that $x,y{\in}\s
N_q\s{Im}\gamma{\subset}\s N_qE$.

A subgroup $H$ of a finitely generated group $G$ is said to be \it
weakly $\alpha$-quasiconvex \rm if there is a proper (i.e.
stabilizers are finite) action of $G$ on a connected graph
$\Gamma$ such that some $H$-orbit ${\subset}\Gamma^0$ is weakly
$\alpha$-quasiconvex.

\medskip

We precise that the word 'weakly' appears in the above definition
since we do not request the above property to be true for every
path having endpoints in $E$ (in which case it is called
$\alpha$-quasiconvex).

 The main result of this Section relates the (weak)
 $\alpha$-quasiconvexity (see 4.1)
   with the existence of
 cocompact action   outside of the limit set  (see 4.2). The constant
 $\lambda_0$ below is fixed in our Convention 2 (see \ref{Flmap}).

\vskip5pt \bf Theorem B\sl.\ Let a finitely generated group $G$
act 3-discontinuously and 2-cocompactly on a compactum $\ti T$.
Let
  $\s{Par}$ be the set of the parabolic points for this action.
  Suppose that
$A{=}\ti T\setminus T{\ne}\varnothing$ where $T{=}\b\Lambda G$ is
the limit set for the action. Then there exists a constant
$\lambda_0\in]1,+\infty[$ such that the following properties of a
subgroup $H$ of $G$ are equivalent:\hfil\penalty-10000 $\s a:$$H$
is weakly $\alpha$-quasiconvex for some
 distortion function $\alpha$ for which $\alpha(n){{\leqslant}}\lambda_0^n\ (n\in\N),$ and for
every $\g p{\in}\s{Par}$ the subgroup $H{\cap}\s{St}_G\g p$ is
either finite or has finite index in $\s{St}_G\g
p$;\hfil\penalty-10000 $\s b:$ the space $(\ti T\setminus
\b\Lambda H)/H$ is compact;\hfil\penalty-10000 $\s c:$ for every
distortion function $\alpha$ bounded by $\lambda_0^n\ (n\in\N),$
every $H$-invariant $H$-finite set $E{\subset}A$ is
$\alpha$-quasiconvex and  for every $\g p{\in}\s{Par}$ the
subgroup $H{\cap}\s{St}_G\g p$ is either finite or has finite
index in $\s{St}_G\g p$\rm.\hfill$\blacklozenge$

\bigskip

Note that the implication `$\s{c{\Rightarrow}a}$' is trivial. The
implication `$\s{b{\Rightarrow}c}$' is rather simple (see
\ref{easyImpl} below). The Section is mainly devoted to the proof
of `$\s{a{\Rightarrow}b}$'.

\subsection{Preliminary results.}

We start with the following  obvious:

\begin{prop}\label{actOnSet}
Let a group $G$ act properly on a set $M$. Let $A_0,A_1$ be
subgroups of $G$ and let $E_\iota$ be $A_\iota$-finite \bf
non-disjoint \sl subsets of $M$ for $\iota{=}0,\mkern-5mu 1$. Then
$|E_0{\cap}E_1|{=}\infty\Leftrightarrow|A_0{\cap}A_1|{=}\infty$.\qed\end{prop}

Since now on we fix a discrete finitely generated group $G$, a
compactum $\ti T$ and a 3-discontinuous  2-cocompact action
$G{\on}\ti T$. Denote by $T{=}\b\Lambda G$ the limit set and
suppose that $A{\leftrightharpoons}\ti T\setminus \b\Lambda
G{\ne}\varnothing$. Since $G$ is finitely generated and $A$ is
$G$-finite there is a $G$-finite set
$\Gamma^1{\subset}\b\Theta^2A$ such that the graph $\Gamma$ with
$\Gamma^0{=}A$ is connected (see e.g. \cite[Lemma 3.11]{GP10}). We
  fix the graph $\Gamma$.

For $x,y{\in}A$ denote by $[x,y]$ the geodesic $\{a{\in}A:\s
d(x,a){+}\s d(a,y){=}\s d(x,y)\}$ between $x$ and $y.$
\begin{prop}\label{dst2pr} Let $E$ be any subset of $A$ and let
 $x{\in}A,z{\in}E,y{\in}\s{Pr}_Ex$ (see Subsection 3.2). Then $\s d(z,[x,y])\geqslant{1\over2}\s
d(z,y)$.\end{prop}\sl Proof\rm. If $t{\in}[x,y]$ and $\s
d(z,t){=}r$ then $\s d(t,y){\leqslant}r$ by the definition of
$\s{Pr}_E$.\hfil\penalty-10000 By $\triangle$-inequality, $\s
d(z,y){\leqslant}\s d(z,t){+}\s d(t,y){\leqslant}2r$.\qed

\medskip

\begin{prop}\label{clPar}
Let $O$ be a weakly $\alpha$-quasiconvex $H$-orbit in $A$. % satisfying condition $\s a$ of \ref{mainThm}.
A parabolic point $\g p{\in}\s{Par}$ belongs to $\overline O$ if
and only if $H{\cap}P$ is infinite where $P{=}\s{St}_G\g
p$.\end{prop}\sl Proof\rm. If $H{\cap}P$ is infinite then
obviously $p\in\overline O$. Suppose $p\in\overline O$. Let us fix
$v{\in}O$. By compactness argument there is an $\alpha$-ray
$\gamma:\Bbb Z_{\geqslant0}\to A$ starting in $\s N_qv$ and
converging to $\g p$ whose image is contained in $\s N_qO$. Let
$S$ be a system of horospheres and let $d{\leftrightharpoons}\s
d(\gamma(0),S_{\g p})$. So $\gamma(0)$ belongs to the $P$-finite
set $\s N_dS_{\g p}$. Hence $\s{Im}\gamma$ is contained in a
$P$-finite set $\s H_\alpha\s N_dS_{\g p}$. Thus the intersection
of $O$ with the $P$-finite set $\s N_q\s H_\alpha\s N_dS_{\g p}$
is infinite. By \ref{actOnSet}, $H{\cap}P$ is infinite. \qed

\medskip

For a fixed Floyd function $f$ we fix an appropriate distortion
function  $\alpha$ (see \ref{Ka}).
\begin{prop}\label{newLemma}
Let $H$ be a subgroup of $G$ acting cocompactly on $\ti T\setminus
\b\Lambda H,$ and let $E$ be an $H$-invariant $H$-finite subset of
$A.$  Then there exist constants $R,d$ such that, for every
$\alpha$-distorted path $\gamma$ in $A$, if the distance between
the $E$-projections of its endpoints is greater than $R$ then
these projections are contained in $\s
N_d\s{Im}\gamma$.\end{prop}\sl Proof\rm. Let $K$ be a compact
fundamental set for the action of $H$ on $\ti T\setminus \b\Lambda
H.$ By  \ref{prOfA} the set $\s{Pr}_EK$ is finite so the number
$\varepsilon_0{\leftrightharpoons}{1\over3}\s{min}\{\overline\d_z(K,\partial
E):z{\in}\s{Pr}_EK\}$ is positive.\hfil\penalty-10000 We take
$\varepsilon{\in}(0,\varepsilon_0)$ and $d{\leftrightharpoons}\s
K_\alpha(\varepsilon)$ where $\s K_\alpha$ is the Karlsson
function  (\ref{KaFn}) corresponding to $\alpha.$ Let
 $R{\leftrightharpoons} \s{max}\{\s
C_E(\varepsilon),2d\}$  where $C_E(\varepsilon)$ is finite as $E$
is weakly homogeneous  (see Subsection 3.3).

Let $\gamma:[0, N]\to A$ be an $\alpha$-distorted path. Up to
applying an element of $H$ we can suppose that  $\gamma(0){\in}K$.
 Let $z{\in}\s{Pr}_E\gamma(0)$, $t{\in}\s{Pr}_E\gamma(N)$.
\begin{center}
\begin{picture}(100,160)(0,-105)
\put(-50,50){\circle*{5}}%x
\put(-63,47){\makebox(0,0)[c]{$\gamma(0)$}}
\put(-12,12){\circle*{5}}%z
\put(-6,8){\makebox(0,0)[c]{$z$}}
\put(40,16){\makebox(0,0)[c]{$E$}}
\put(-41,-13){\makebox(0,0)[c]{$\gamma$}}
\put(105,14){\makebox(0,0)[c]{$\b\Lambda G$}}
\qbezier(-24,-25)(-24,24)(25,24)
\qbezier(-24,-25)(-24,-74)(-3,-100)
\qbezier(95.5,-10)(74,24)(25,24)
\put(-22,-51.5){\circle*{5}}%t
\put(-99,-90){\line(2,1){78}} \put(-15,-49){\makebox(0,0)[c]{$t$}}
\put(-99,-90){\circle*{5}}%y
\put(-116,-90){\makebox(0,0)[c]{$\gamma(N)$}}
\put(-50,50){\line(1,-1){38}} \qbezier(-50,50)(-20,10)(-30,0)
\qbezier(-30,0)(-40,-10)(-35,-30)
\qbezier(-99,-90)(-30,-50)(-35,-30) \thicklines
\qbezier(-8,-102)(77,-77)(102,8)%\Lambda G
\end{picture}\end{center}

Since $d(z,t) > C_E(\varepsilon)$  so by  \ref{CEa} and \ref{CE}
we obtain $\overline\d_z(t,\partial E){\leqslant}\varepsilon$. We
have $d(z, t)\geq R\geq 2d$ so
    $d(z, [t, \gamma(N)])\geq d$ by Lemma \ref{dst2pr}. Thus
$\overline\d_z(\gamma(N),t){\leqslant}\varepsilon$. If by
contradiction $d(z,\s{Im}\gamma) > d$ then  we would have
$\overline\d_z(\gamma(0), \gamma(N)){\leqslant}\varepsilon.$ So
summing all these three inequalities we would have
 $\overline\d_z(K,\partial
E){\leqslant} \overline\d_z(\gamma(0), \partial
E){\leqslant}\overline\d_z(\gamma(N),\gamma(0)){+}\overline\d_z(t,\gamma(N)){+}
\overline\d_z(t,\partial E) {\leqslant}
3\varepsilon{<}3\varepsilon_0$ contradicting to the choice of
$\varepsilon_0$.\qed \vskip3pt

As an immediate corollary we obtain the following

\begin{prop}\label{oldLemma}
Given a system $S$ of horospheres there exist constants $R,d$ such
that for every $\g p{\in}\s{Par}$ and every $\alpha$-distorted
path $\gamma$ in $A$ if the distance between the $S_{\g
p}$-projections of its endpoints is greater than $R$ then these
projections are contained in $\s N_d\s{Im}\gamma$.\end{prop}\sl
Proof\rm. Since $\s{Par}$ is $G$-finite the problem reduces to the
case of  a fixed $\g p{\in}\s{Par}$. The subgroup
$P{\leftrightharpoons}\s{St}_G\g p$ acts cocompactly on $\ti
T\setminus \g p$. So the assertion follows from \ref{newLemma}
applied to $H{=}\s{St}_G\g p$ and $E{=}S_{\g p}$.\qed
\subsection{Implication `$\s{b{\Rightarrow}c}$'.}\label{easyImpl} We fix a
subgroup $H$ satisfying property `$\s b$' of Theorem B. Let $E$ be
an $H$-finite $H$-invariant subset $E$ of $A$. Denote by
$\overline E$ the closure of $E$ in $\ti T$ and $\partial
E{\leftrightharpoons}\overline E\setminus E$. It follows from the
the convergence property  that $\partial E{=}\b\Lambda H$. The
$\alpha$-quasiconvexity of $E$ follows from Proposition
\ref{subgrConvex}. So the only thing to prove is the last part of
the statement.

We fix a compact fundamental domain $K$ for the action $H{\on}(\ti
T\setminus \partial E)$. Consider  a system of horospheres
$\{S_{\g p} : {\g p}{\in}\s{Par}\}$ (see  5.1). The set $\Cal
P_H{=}\{\g p{\in}\s{Par}:|S_{\g p}\cap E|{=}\infty\}$ is
$H$-invariant. As $K{\cap}\partial E{=}\emptyset$ there exists an
entourage   $\b u\in\Theta^2\ti T$  such that $\b u{\cap}
(K{\times}\partial E){=}\varnothing$. Since every parabolic
subgroup $\s{St}_G\g p$ is quasiconvex (see \ref{subgrConvex}), it
is dynamically quasiconvex (Theorem A). Furthermore there are at
most finitely many $G$-non-equivalent parabolic points \cite[Main
Theorem, a]{Ge09}. Hence the set $\Cal P_{\b
u}{\leftrightharpoons}\{\g p{\in}\s{Par}:S_{\g p}$ is not $\b
u$-small\} is finite. We have $K\cap{\cup}\{S_{\g p}:\g p{\in}\Cal
P_H\}{=}K\cap{\cup}\{S_{\g p}:\g p{\in}\Cal P_{\b u}{\cap}\Cal
P_H\}$. For every $\g p{\in}\ \Cal P_H$ the set $K{\cap}S_{\g p}$
is finite as otherwise it would contain the unique limit point $\g
p$ of the infinite set $S_{\g p}\cap E$ which is impossible. So it
follows that the set ${\cup}\{S_{\g p}:\g p{\in}\Cal P_H\}$ is
$H$-finite. By \ref{horosph}.$\s b$ each element of $A$ belongs to
at most finitely many horospheres, so the set ${\Cal
S}_H{=}{\cup}\{\g p{\times}S_{\g p}: \g p{\in}\Cal P_H\}$ is
$H$-finite too.

By the property \ref{horosph}.$\s c$, applied to ${\Cal S}_H$,
each $S_{\g p}$
 is $H{\cap}\s{St}_G\g p$-finite ($\g p{\in}\Cal P_H$). It implies
  that the index of $H{\cap}\s{St}_G\g p$ in $\s{St}_G\g p$
is finite.

\subsection{Implication `$\s{a{\Rightarrow}b}$'}  We fix a subgroup $H$ of
$G$ satisfying condition `$\s a$' of Theorem B. Let $\Gamma$ be a
locally connected graph where $G$ acts properly. Let $O$ be a
weakly quasiconvex $H$-orbit satisfying   `$\s a$'  with the
parameter $q$. Denote by $A$ the vertex set $\Gamma^0$ and by
$\widetilde T$ the union $A{\sqcup}\b\Lambda G$.

\begin{prop}\label{hFinite}
For every system of horospheres $S$ there exists a constant $c$
such that $S_{\g p}{\subset}\s N_cO$ for every $\g p{\in}\overline
O{\cap}\s{Par}$.\end{prop}

\sl Proof\rm. Since the subgroup $H$ is weakly $\al$-quasiconvex
in $G$ there exists an $\alpha$-isometric map $\varphi {:} H{\to}
\Gamma$ with the distortion function $\alpha$. Then by
\cite[Theorem C]{GP09} the subgroup $H$ is relatively hyperbolic
with respect to the system $\{H{\cap}\s{St}_G\g p:\g
p{\in}\s{Par}\}$. Thus this system of maximal parabolic subgroups
of $H$ contains at most finitely many $H$-conjugacy classes
\cite[Main Theorem, a]{Ge09}.

By Proposition \ref{clPar} for every $\g p{\in}\overline
O{\cap}\s{Par}$ the set $P{\cap}H$ is infinite where
$P{=}\s{St}_G\g p$. Then by our assumption    $\vert
P{:}H{\cap}P\vert {<}\infty$. Since there are at most finitely
many distinct $H$-conjugacy classes of such subgroups the set of
all indices $\vert \s{St}_G\g p : H{\cap}\s{St}_G\g p\vert\ (\g
p{\in}\s{Par}) $ is
  bounded. The Proposition follows. \qed
\medskip

 \bf Remark.\rm\ In the above proof  the fact that the subgroup $H$ is
 relatively hyperbolic itself and can contain at most finitely
 many conjugacy classes of distinct parabolic subgroups was
 essential. In general there are examples of geometrically
 finite Kleinian groups containing finitely generated subgroups
 having infinitely many conjugacy classes of parabolic subgroups
 \cite{KP91}.

\medskip

\bf Definition\rm. Let $O\subset A$ be an $H$-orbit of a point
$v{\in}A.$ The set
\setcounter{equation}1%
 \begin{equation}\label{diric}F_v{\leftrightharpoons}\{x{\in}A :\s d(x,O){=}\s
 d(x,v)\}\end{equation}
 is called Dirichlet set at   $v.$
\medskip

 \bf Remark.\rm\ The set $F_v$ is a discrete analog of the
 Dirichlet fundamental set for a discrete subgroup of the isometry
 group of the real hyperbolic space.

\setcounter{prop}2%
\begin{prop}\label{diricfund} The set $F_v$ is a $v$-star convex
fundamental set for the action of $H$ on $A.$\end{prop}

\sl Proof\rm. For every point $x{\in} A$ there exists $w{=}h(v)\in
O$ such that $d(x, w){=}d(x, O)$. So $h^{-1}x{\in} F_v$ and
$A{=}\bigcup_{h\in H}hF_v.$

We have $v{\in} F_v$. To show that $F_v$ is $v$-star convex we
need to show that if $w{\in} F_v$ then for any $t{\in} [w,v]$ we
have $t{\in} F_v.$  Suppose not then there exists $u{\in}
O\setminus  v$ such that $d(t,u){ < }d(t,v)$. Then
$d(w,v){=}d(t,v)+d(t,w){>} d(t,u) + d(t,w){\geq} d(u,w)$ which is
impossible as $w{\in} F_v.$\qed

\medskip

  The main step in proving   the implication
`$\s{a{\Rightarrow}b}$' is the following.

\begin{prop}\label{fundDom} The closure of the set $F_v$ is
disjoint from $\partial O$.\end{prop}

\vskip3pt

\bf Corollary.\sl \ The closure $\overline F_v$ of $F_v$ in $\ti
T$ is a compact fundamental set for the action of $H$ on $\ti
T\setminus{\b\Lambda H}.$

\vskip3pt

\sl Proof of Corollary\rm.  Let $x\in\ti T\setminus{\b\Lambda H},$
we need to show that there exists $h\in H\ :\ h(x)\in \overline
F_v.$ If $x\in A$ then it follows from \ref{diricfund}.

Let $x\in {\b\Lambda G}\setminus {\b\Lambda H}.$ Then there exists
a sequence $(x_n)\subset A$ converging to $x.$ Let $h_n\in H$ such
that $h_n(x_n)\in F_v.$  Suppose first that the set
$\{h^{-1}_n(F_v)\}_n$ is infinite. By \ref{fundDom} $\overline
F_v{\cap}{\b\Lambda H}{=}\emptyset$ so up to passing to a
subsequence we obtain $h_n^{-1}(y)\to x$ for every $y\in F_v.$
Then $x\in {\b\Lambda H}$ which is impossible. So the
  set $\{h^{-1}_n(F_v)\}_n$ is finite and up to a new
  subsequence we have $x_n\in h^{-1}F_v$ for
  a fixed
  $h\in H.$ Thus   $h(x)\in \overline F_v.$\qed

\medskip

\noindent To prove \ref{fundDom} we need the following.
\setcounter{lemma}4
\begin{lemma}\label{prHoro}Let $S$ be a system of horospheres for the action $G{\on}\ti T$.
There exists a constant $D_S$ such that, for $\g p{\in}\overline
O{\cap}\s{Par}$, $\s{diam_d\,Pr}_{S_{\g
p}}F_v{\leqslant}D_S$.\end{lemma}

\sl Proof\rm. Let $w{\in}F_v$. Suppose that the distance between
the projections of $v$ and $w$ is greater than the constant  $R$
of \ref{oldLemma}. \vskip2pt
\begin{center}\begin{picture}(120,105)(0,-40)
\put(-60,0){\line(1,0){120}} \put(-40,40){\line(0,-1){40}}
\put(-40,40){\circle*{5}}%v
\put(-40,44){\makebox(0,0)[b]{$v$}} \put(40,60){\line(0,-1){90}}
\put(40,60){\circle*{5}}%w
\put(40,64){\makebox(0,0)[b]{$w$}} \qbezier(-40,40)(30,-20)(40,60)
\put(-40,0){\circle*{5}}%v
\put(-40,-2){\makebox(0,0)[t]{$\s{Pr}_{S_{\g p}}v$}}
\put(40,0){\circle*{5}}%v
\put(43,-2){\makebox(0,0)[lt]{$\s{Pr}_{S_{\g p}}w$}}
\put(40,-30){\circle*{5}}%v
\put(44,-30){\makebox(0,0)[l]{$x{\in}O$}}
\put(22,18){\circle*{5}}%v
\put(40,0){\line(-1,1){18}} \put(-10,20){\makebox(0,0)[b]{$m$}}
\put(39,-15){\makebox(0,0)[c]{$c$}}
\put(30,10){\makebox(0,0)[c]{$d$}}
\put(19,24){\makebox(0,0)[c]{$t$}}
\put(0,-19){\makebox(0,0)[c]{$S_{\g p}$}}
\end{picture}\end{center}
Then there exists $d$ such that $\s d(\s{Pr}_{S_{\g
p}}w,[w,v]){=}\s d(\s{Pr}_{S_{\g p}}w,t){\leqslant}d$ for some
$t{\in}[w,v]$. Let $m{\leftrightharpoons}\s d(t,v)$ and let $c$ be
the constant  from Proposition \ref{hFinite}. Then by Proposition
\ref{diricfund} $t{\in}F_v$, so for every $x\in O$  we have
\hfil\penalty-10000 $m{\leqslant}d(t,x){\leqslant}d{+}c$, $\s
d(v,S_{\g p}){\leqslant}m{+}d{\leqslant}2d{+}c$; $\s{d(Pr}_{S_{\g
p}}v,\s{Pr}_{S_{\g p}}w){\leqslant}d(\s{Pr}_{S_{\g p}}v,v){+}d(v,
\s{Pr}_{S_{\g p}}w))
\\ {\leqslant}2d{+}c{+}m{+}d{\leqslant}4d{+}2c{\leftrightharpoons}D_S$.\qed
\vskip3pt \sl

 \vskip3pt
Proof \rm of \ref{fundDom}. Suppose that the assertion is false,
and let  $\g t{\in}\overline{F_v}{\cap}\partial O$.

By compactness argument there exists an infinite $\alpha$-geodesic
$\gamma$ starting at $v\in O$ and converging to $\g t$. It is the
limit of a sequence of $\alpha$-geodesics whose endpoints are $v$
and $\g t_n{\in} O$ such that $\g t_n\to\g t.$ So by the
$\alpha$-weak quasiconvexity we can assume  that $\s{Im}\gamma$ is
contained in $\s N_qO.$

Choose $v_n{\in}\s{Im}\gamma$   such that $\s d(v_n,v){>}2n{+}q$.
We claim that $\s d(v_n,F_v){>}n$. Indeed let $o_n\in O$ be such
that $d(o_n, v_n){\leqslant}q$. Take $v'_n\in \s{Pr}_{F_v}v_n.$
Then

$d(v'_n, v){\leqslant} d(v'_n, o_n){\leqslant} d(v'_n,
v_n){+}d(v_n,o_n){\leqslant}d(v'_n, v_n)+q.$

Hence
$2n+q{<}d(v,v_n){\leqslant}d(v,v'_n)+d(v'_n,v_n){\leqslant}2d(v'_n,v_n)+q$
and the claim follows.

 Since $\s N_qO$ is $H$-finite there exist $h_n{\in}H$ such that $h_nv_n$
belong to a ball of a finite radius  centered at $v$. By passing
to a subsequence we can reduce the situation to the case when
$h_nv_n{=}w$ independently on $n$. Moreover we can assume that the
sequence $h_n\gamma$ converges in the Tikhonoff topology to an
infinite $\alpha$-geodesic $\beta:\Bbb Z\to A$ with
$\beta(0){=}w$.

We have $d(w, h_nF_v){>}n$. By Proposition \ref{diricfund}
$h_nF_v$ is $h_nv$-star convex. Then by Proposition \ref{gKaFn} we
obtain that $\s{diam}_{\overline\d_w}(h_nF_v)\to 0$ where
$\overline\d_w$ is the shortcut metric with respect to a Floyd
function forming with $\alpha$ an appropriate pair. Thus
$\overline\d_w(h_nv,h_n\g t)\to 0$ and  $\beta$ is an
$\alpha$-horocycle based at a parabolic point $\g q{\in}\s{Par}$
(see \ref{horoc}). Then by \ref{horosph}
 $\s{Im}\beta$ is contained in the $Q$-finite set $S_{\g q}{=}\s H_\alpha\g
q$ where  $Q{\leftrightharpoons}\s{St}_G\g q$. It follows that the
intersection of $O$ with the $Q$-finite set $\s N_q\s H_\alpha\g
q$ is infinite. By \ref{actOnSet} the subgroup $H{\cap}Q$ is
infinite and so $\g q{\in}\overline O.$

By   Lemma  \ref{prHoro}     the diameter of $\s{Pr}_{S_{\g
q}}(h_nF_v){=}\s{Pr}_{h_n^{-1}S_{\g q}}F_v$ is bounded. Since the
diameter of $(h_n\s{Im\gamma){\cap}Im}\beta$ tends to infinity
 there exists a point
$w_N{\in}(h_n\s{Im\gamma){\cap}Im}\beta$ such that the number
$N{=}d(w_N, \s{Pr}_{S_{\g q}}(h_nF_v))$ is arbitrarily large. To
obtain a contradiction we choose $N$  in few steps. The endpoints
of the curve $h_n\gamma$ belong to $h_n\overline F_v$, so we have
$\s{Pr}_{S_{\g q}}(\partial h_n\s{Im\gamma})\subset\s{Pr}_{S_{\g
q}}(h_nF_v).$ Since $w_N{\in} S_{\g q}$ we first assume  that $N
> R$ where $R$ is the constant from Proposition \ref{oldLemma}.
Then it implies that there exists $d{=}d(R)$ such that both
$\alpha$-distorted subpaths of $h_n\gamma$ joining $w_N$ with its
endpoints meet the $d$-neighborhood $U{=}\s N_d(\s{Pr}_{S_{\g
q}}(h_nF_v))$ of $\s{Pr}_{S_{\g q}}(h_nF_v)$ in $A.$ On the
picture below the points $w$ and $t$ belong to $\s{Pr}_{S_{\g
q}}(h_nF_v)$ and are close to the subpaths of $h_n\gamma.$
\begin{center}\begin{picture}(100,145)(0,-50)
\renewcommand{\qbeziermax}{1000}
\qbezier(40,0)(40,16.6)(28.3,28.3)
\qbezier(0,40)(16.6,40)(28.3,28.3)
\qbezier(-40,0)(-40,16.6)(-28.3,28.3)
\qbezier(0,40)(-16.6,40)(-28.3,28.3)
\qbezier(-40,0)(-40,-16.6)(-28.3,-28.3)
\qbezier(0,-40)(-16.6,-40)(-28.3,-28.3)
\qbezier(40,0)(40,-16.6)(28.3,-28.3)
\qbezier(0,-40)(16.6,-40)(28.3,-28.3)
\put(0,40){\circle*{5}}%w
\put(0,44){\makebox(0,0)[b]{$w_N$}}
\put(0,-40){\circle*{5}}%\g q
\put(0,-48){\makebox(0,0)[c]{$\g q$}}
\put(-80,0){\circle*{5}}%p
\put(-83,0){\makebox(0,0)[r]{$h_nv$}}
\put(-40,0){\circle*{5}}%p
\put(-36,0){\makebox(0,0)[l]{$w$}}
%\put(-37.5,14){\circle*{5}}%y_-
%\put(-34,12){\makebox(0,0)[l]{$y_-$}}
\put(18,84){\makebox(0,0)[t]{$h_nF_v$}}
\put(4,33){\makebox(0,0)[tl]{$\s H_\alpha\g q$}}
%\put(-26,23){\makebox(0,0)[l]{$ht_k$}}
\put(-52,10){\makebox(0,0)[rb]{$h_n\gamma$}}
\put(-80,0){\line(3,1){42.5}} \qbezier(-80,0)(-80,30)(-60,60)
%\put(-60,60){\circle*{5}}%x_k
%\put(-65,66){\makebox(0,0)[b]{$hx_k$}}
%\put(-60,60){\line(1,-1){31.7}}
%\put(-28.3,28.3){\circle*{5}}%t
\qbezier(-60,60)(-40,90)(-10,90) \qbezier(-10,90)(40,90)(65,65)
\qbezier(80,20)(80,50)(65,65) \qbezier(80,20)(80,0)(49,-20)
\put(49,-20){\circle*{5}}%p
\put(52,-22){\makebox(0,0)[tr]{$h_n\g t$}}
\put(38.6,-10){\circle*{5}}%p
\put(36,-9){\makebox(0,0)[r]{$t$}}
%\put(39,9){\circle*{5}}%y_+
\put(45,-7){\makebox(0,0)[lb]{$h_n\gamma$}}
\put(49,-20){\line(-1,3){12}}
\qbezier(-45,1)(-45,10)(-20,7) \qbezier(-45,1)(-45,-10)(-20,-13)
\qbezier(20,1)(5,4)(-20,7) \qbezier(20,1)(45,-4.3)(43,-10)
\qbezier(20,-15)(5,-16)(-20,-13) \qbezier(20,-15)(43,-18)(43,-10)
\put(0,-7){\makebox(0,0)[c]{$\s{Pr}_{S_{\g q}}(h_nF_v)$}}
\put(70,-24){\makebox(0,0)[c]{$\b\Lambda G$}} \thicklines
\qbezier(-70,-30)(0,-50.3)(70,-30)%\Lambda G

\end{picture}
\end{center}

%\hspace*{2.5cm}{\rm Bounded projection on the horosphere $S_{\g
%q}{=}\s H_\alpha\g q$. $w\in\s{Pr}_{S_{\g q}}h_nv,\ \s{Pr}_{S_{\g
%q}}h_n\gamma.$

%\bigskip

Therefore $h_n\s{Im\gamma}\subset \s H_\alpha (U)$. By Proposition
\ref{hullLoc} there exists $s{=}s(d,\alpha)$ such that
$h_n\s{Im\gamma}\subset \s N_s(\s{Pr}_{S_{\g q}}(h_nF_v))$.
Assuming finally that $N > \s{max}\{R, s\}$ we obtain
$w_N{\in}h_n\s{Im\gamma}\setminus N_s(\s{Pr}_{S_{\g q}}(h_nF_v))$
which is a contradiction. \qed

\medskip

 The condition $`\s{b}$' of Theorem B does not depend on the choice of $A$. We
formulate this independence as follows

\medskip

 \bf Corollary.\sl\   Let a group $G$ act on compacta $X$ and $\ti T$
 3-discontinuously and 2-cocompactly such that the corresponding limit sets
 are
 proper subsets of them.
Let $\varphi:X\to \ti T$ be a continuous equivariant map bijective
on the limit sets. Then a subgroup $H{<}G$ acts cocompactly on
$\ti T\setminus \b\Lambda H$ if and only if $H$ acts cocompactly
on $X\setminus \b\Lambda H$.\qed

\section{Subgroups of convergence groups with
 proper limit sets.}
\subsection{Dynamical boundness of subgroups.}
Let $G$ be a group acting 3-discontinuosly on a compactum $\ti
T{=}T{\sqcup}A$ where $T{=}\b\Lambda G$  and $A$ is a non-empty,
discrete and $G$-finite set (see Subsection 2.4).

\medskip

 \bf Remark.\rm\ We do not assume in this Subsection that the action is 2-cocompact
nor that $G$ is finitely generated.

\medskip

 \bf Definition.\rm  \  {\it Let $G$  acts 3-discontinuously on a
compactum $T$. A subgroup $H$ of $G$ is called {\it dynamically
bounded} if every infinite set of  elements $S{\subset} G$
contains an infinite subset $S_0$ such that $\displaystyle
T\setminus \bigcup_{s{\in} S_0} s(\b\Lambda H)$ has a non-empty
interior.}

\medskip

\noindent We start by giving several equivalent reformulations of
this notion.

\begin{prop}\label{equivbound}. Let $T$ be a metrisable compactum
and  $G$ be a group  acting 3-discontinuously on $X.$ Then the
following statements are equivalent:

\medskip

\begin{itemize}

\item[\sf 1)] $H$ is dynamically bounded in $G$.

\medskip

\item[\sf 2)] There exist finitely many proper closed subsets
$F_1, ..., F_k$ of $T$ such that $$\forall g{\in} G\ \exists\
i{\in}\{1,...,k\}\ :\ g(\b\Lambda H){\subset} F_i.$$

\medskip

\item[\sf 3)] In the space  ${\rm Cl}(T)$ of closed subsets of a
compactum $T$ equipped with the Hausdorff topology one has
$$T{\not\in} \overline{\{g(\b\Lambda H)\ :\ g{\in} G\}}.$$ \end{itemize}

\end{prop}

\bigskip

 \bf Corollary.\sl\   The dynamical boundness is a hereditary
property with respect to subgroups, i.e. if a subgroup $H$ of $G$
is dynamically bounded then any subgroup of $H$ is so.

\medskip

\sl Proof of Corollary\rm.\ It follows immediately e.g. from the
condition 2). Indeed if $H_0< H$ then $\La H_0{\subset} \La H$ and
so the sets $F_i$ existing for $H$ work equally for $H_0\
(i{=}1,...,k).$ \qed

\medskip

\sl Proof of the Proposition.\rm\ Let us prove the following
implications : $2)\Rightarrow 1)\Rightarrow 3)\Rightarrow 2).$

\medskip $2)\Rightarrow 1).$ Let $S{\in} G\setminus H$ an
infinite set of pairwise distinct elements modulo $H.$ Then there
exists an infinite subset $S_0{\subset} S$  such that $\exists\
i{\in} \{1,...,k\}\ \forall s{\in} S_0\ s(\b\Lambda H){\subset}
F_i$. The set $F'_i{=}T\setminus F_i$ is open so we are done.

\medskip

$1)\Rightarrow 3).$ Suppose by contradiction that 3) is not true
and there exists a sequence $(g_n){\subset} G$ such that
$g_n(\b\Lambda H)\to T$ in the Hausdorff topology. Then the same
is true for any its subsequence contradicting the condition 1).

\medskip

$3)\Rightarrow 2).$ We provide a topological proof. Recall first
that every entourage $\bu{\in}\ent$ defines a distance function
$\Delta_{\bu}$ on $T$ which is the maximal one with the property
$(x,y){\in}\bu{\cap}\Theta^2T$ if and only if
$\Delta_{\bu}(x,y){\leqslant}1$ (see e.g. \cite{GP10}). So for
every $\bu{\in}\ent$ we define the entourage $\bw{=}\bu^k\
(k{\in}\N)$ such that $(x,y){\in}\bw$ if and only if
$\Delta_{\bu}(x,y){\leqslant}k.$

 Let now
$\bu{\in}\ent$ be an entourage on $T.$ For any subset $C{\subset}
T$  its $\bu$-neighborhood $C\bu$ is the set $\{x{\in} T\ \vert\
\exists\ y{\in} C : (x,y){\in}\bu\}.$ By the condition 3) there
exists $\bu{\in} \ent$ such that $\forall g{\in} G\ g(\b\Lambda
H)\bu{\not=}T.$ In other words $\forall g{\in} G\ \exists \
p_g{\in} T : \forall\ y{{\in}} g(\b\Lambda H)\ (p_g,y){\not\in}
\bu$ (i.e. $p_g\bu\cap g(\b\Lambda H){=}\emptyset).$ Take an
entourage $\bv{\in}\ent$ such that $\bv^2{\subset} \bu$ meaning
that $(x,y){\in}\bv$ and $(y,z){\in} \bv$ implies $(x,z){\in}
\bu.$ Since $T$ is compact there exists a finite $\bv$-net $\cal
P{\subset} T$ such that $\forall x{\in} T\ \exists\ y{\in} \cal P\
:\ (x,y){\in}\bv.$ So for every $g\in G$ there is $q_g{\in}\cal P\
:\ (p_g, q_g){\in} \bv.$ It follows that $q_g\bv\cap g(\b\Lambda
H){=}\emptyset$ as otherwise $p_g\bu\cap g(\b\Lambda
H){\not=}\emptyset$. The set $F_q{=}(q_g\bv)'$ is the desired
closed subset of $T.$ \qed

\medskip

\bf Remark.\rm\ In the above proof we need the metrisability of
$T$ only to prove the second implication as the choice of a
sequence converging to an accumulation point in a topological
space without countable basis is not possible in general.

\medskip

\begin{prop}\label{dcb} If $H<G$ is dynamically quasiconvex then
it is dynamically bounded. \end{prop}

\bf Proof.\rm\ Let us fix an entourage $\bu{\in} \ent$ such that
$T$ is not $\bu^4$-small (i.e. $\displaystyle {\rm
diam}_{\Delta_{\bu^4}}(T){>1}$). Then by compactness of $T$ there
exists a finite $\bu$-net $\cal P.$ So for any $x{\in} T\ \exists\
y{\in} \cal P$ such that $(x,y){\in}\bu.$ Let $S{\subset} G$ be an
infinite set of elements. Then there is an infinite subset
$S_0{\subset} S$ such that $\exists y{\in}\cal P\ \forall s{\in}
S_0\ y{\in} s(\b\Lambda H)\bu$. Since $H$ is dynamically
quasiconvex up to removing a finitely many elements we can assume
that for all $s{\in}S_0$ we have $s(\b\Lambda H)$ is $\bu$-small.
Therefore $\forall\ s{\in} S_0\ :\ s(\b\Lambda H){\subset} U_y$
where $U_y$ is an $\bu^2$-small neighborhood of $y.$

Then there exists $z{\in} \cal P\setminus y$ having an
$\bu^2$-small neighborhood $U_z$ such that $U_y\cap
U_z{=}\emptyset.$ Indeed otherwise every point of $T$ would belong
to an $\bu^4$-small neighborhood of $y$ which is impossible. So
$T\setminus \bigcup_{s{\in} S_0}s(\b\Lambda H)$ has a non-empty
interior. \qed

\medskip

 The following Proposition shows that  a
dynamically bounded subgroup acting cocompactly outside  the limit
set on $T$ do the same on $\ti T$.

\begin{prop}\label{extens} Let  $G$ act 3-discontinuously on
$\widetilde T{=}T{\sqcup}A$. Suppose $H$ is a dynamically bounded
subgroup of $G$ acting cocompactly on $T\setminus \b\Lambda H$.
Then   $H$ acts cocompactly on $\widetilde T\setminus \b\Lambda
H.$

\end{prop}

\medskip

By \ref{dcb} every dynamically quasiconvex subgroup is dynamically
bounded so we have.

\medskip

 \bf Corollary.\sl\ Let $T$ and $G$ be as above.
  Let $H< G$ be a dynamically quasiconvex
  subgroup of $G$ acting cocompactly on $T\setminus\b\Lambda H$
  then $H$ acts cocompactly on $\ti T\setminus\b\Lambda H.$ In
  particular if $H$ is a parabolic subgroup for the action of $G$
  on $T$ then it is so for the action on $\ti T$.

\medskip

  \bf Remark\rm.\ If one assumes in addition
   that the action $G\act T$ is 2-cocompact
  then the latter fact also follows from
 \cite[Corollary, 7.2]{Ge09}.

 \medskip

\it Proof of the Proposition.\rm\ Suppose this is not true. Since
$(T\setminus\b\Lambda H)/H$ is compact there exists $H$-invariant
subset $W$ of $A$ such that $\vert W/H\vert{=}\infty$ and all
limit points of $W$ are in $\b\Lambda H.$ The set $A$ is
$G$-finite, so we can assume that $W$ is an orbit $Sa\ (a\in A)$
where $S$ is an infinite set of elements of $G$ representing
distinct right cosets   $H\backslash G$. Since $H$ is dynamically
bounded, $S$ admits an infinite subset $S_0$ such that $
C{=}T\setminus\bigcup_{g{\in} S_0} g^{-1}(\b\Lambda H)$ has a
non-empty interior. Choose $x{\in} C$ which admits a neighborhood
$U_x{\subset}C.$

For every $g{\in} S_0$ we have $g(x){\not\in} \b\Lambda H$, so
there exists $h_g{\in} H$ that $ h_g(g(x)){\in} K,$ where $K$ is a
compact fundamental set for the action $H\act(T\setminus\b\Lambda
H)$. The set $S_1{=}\{\ga_g : \ga_g{=}h_gg,\ g{\in} S_0\}$ is
infinite, so it admits a limit cross $(r, a)^\times{=}r\times
T\cup T\times a$ where $r$ and $a$ are respectively repeller and
attractor points \cite{Ge09}. By our assumption we have $a{\in}
\b\Lambda H.$

 We now claim that $r{\not=}x$.
Suppose not. If  first  there exists $b{\in}\b\Lambda H\setminus
\{a\}$ then we can find $\ga_g{\in} S_1$ close to $(r,a)^\times$
such that $\gamma_g^{-1}(b){\in} U_x.$ By the choice of $U_x$ it
is impossible. So we must have  $\b\Lambda H{=}\{a\}$. Then $a$ is
a parabolic point for the convergence action of $G\act T$
\cite[Proposition 3.2]{Bo99}, \cite[Theorem 3.A]{Tu98} (we note
that the argument of these papers
 can be
applied without assuming the metrisability of $T$). From the other
hand we have $\ga_g^{-1}(a){\not\in} U_x$ and $\ga_g^{-1}(b){\in}
U_x$ for all $b{\not=} a$ and for all elements $\ga_g{\in} S_1$
close to $(r,a)^\times$ (for which $\ga_g^{-1}$ is close to
$(a,r)^\times$). It follows that $a$ is a conical point for the
action $G\act T$ \cite{Bo99}. This is a contradiction. We have
proved $r{\not=}x.$

For any neighborhood $U_a{\subset}T$ of $a$ we have $\ga_g(x){\in}
U_a$ for some $\ga_g{\in} S_1$. Since $\forall g{\in} S_0\
\ga_g(x){\in} K$ we obtain  $K{\cap}U_a{\not=}\emptyset.$ This is
impossible as $a{\in}\b\Lambda H$ and $K$ is compact in
$T\setminus\b\Lambda H$.\qed

\subsection{Finite presentedness of dynamically bounded subgroups.}

The property to act cocompactly outside the limit set for a
subgroup of a RHG has several consequences which have been
established in Sections 4 and 8. The following Proposition gives
one more property of such subgroups.

\medskip
\begin{prop} \label{finpres} Let $G$ be a group acting
3-discontinuously and 2-cocompactly on a compactum $\ti T.$
Suppose $H$ is a subgroup of $G$ acting cocompactly on $\ti
T\setminus \bold\Lambda H$. If $G$ is finitely presented then $H$
also is.\end{prop}  \it Proof\rm. Let $\Gamma$ be a connected
graph on which $G$ acts discontinuously and cocompactly. It is
rather well-known that there exists a simply connected
2-dimensional $\s{CW}$-complex $C(\Gamma)$ such that
$C(\Gamma)^1{=}\Gamma$ and $G$ acts cocompactly on $C(\Gamma)$.
Since the action of $G$ on $\Gamma$ is not necessarily free we
provide for the sake of completeness a short proof of it. Consider
the Cayley graph $\Gamma_0{=}{\rm Cay}(G, S)$ of $G$ corresponding
to the finite generating set $S$ with finitely many defining
relations. Let $R_1$ be a  maximal subset of $G$-nonequivalent
loops in $\Gamma_0$  corresponding to the elements of $S$. Since
the action of $G$ on both graphs $\Gamma_0$ and $\Gamma$ is
cocompact there is an equivariant finite-to-one quasi-isometry
$\varphi : \Gamma_0\to \Gamma$ which is injective everywhere
outside   the set of the preimages of the fixed points for the
action $G\act\Gamma$.
 Let
$C(\Gamma_0)$ be the Cayley 2-dimensional simply connected
$\s{CW}$-complex obtained by gluing 2-cells to the $G$-orbit of
$R_1.$ Denote by $R_2{\subset} \Gamma_0^0$ a maximal subset of
$G$-non-equivalent points on which $\varphi$ is not injective. We
now construct the 2-dimensional $\s{CW}$-complex $C(\Gamma)$ by
attaching  2-cells to the $G$-orbits of the loops $\varphi(a){\in}
\Gamma^1,$ where $a{\in} R_1$ or $a$ is a path connecting a pair
of points in $R_2$  mapped to the same point of $\Gamma^0.$ The
map $\varphi$ extends   continuously and equivariantly to a
surjective map  between the 2-skeletons $C^2(\Gamma_0)\to
C^2(\Gamma).$ Every loop $\gamma{\in} \Gamma$ is a product
$\displaystyle\prod_i \varphi(a_i)$ where each $\varphi(a_i)$ is
trivial in $C(\Gamma).$ Therefore $C(\Gamma)$ is simply connected
and satisfies the claim above.

We will now construct an $H$-invariant 2-dimensional simply
connected $\s{CW}$-complex $\mathcal E$   such that $\mathcal E/H$
is compact. Let $E$ be an $H$-finite and $H$-invariant subset of
$\Gamma^0$. Set $\mathcal E^0{=} E$. Join by an edge each pair of
vertices of $\mathcal E^0$ situated within a distance at most $C$
where $C$ is the constant from Proposition \ref{qIso}. Denote by
${\mathcal E}^1$ the obtained graph. Let $n$ be the maximal length
of the boundary curves of the 2-cells of $C(\Gamma)$ corresponding
to a finite set of generating relations of $G.$ Attaching now a
2-cell to every closed curve of $\mathcal E^1$ of length at most
$n$ denote by $\mathcal E$ the obtained complex.

Let $\s{pr}_E:\Gamma^0\to E$ denote a single valued branch of the
multivalued map $\s{Pr}_E$ obtained by choosing one element from
the image of each vertex. The map $\s{pr}_E$ extends to a
continuous map $C(\Gamma)^1\to \mathcal E^1$ which sends  the
edges of $\Gamma$ to edges of $\mathcal E$ by \ref{qIso}. The
projection of a path in $\Gamma$ is a path in $\mathcal E^1$. The
map $\s{pr}_E$ is surjective so the graph $\mathcal E^1$ is
connected.

Every 2-cell of $C(\Gamma)$ is bounded by a curve which is the
product of curves  of length at most $n$. By construction its
projection  to $\mathcal E$ is also a trivial loop with the same
property. So the map $\s{pr}_E$ extends to a map
 between the 2-skeletons $C(\Gamma)^2\to \mathcal E^2$.

 The complex
 $\mathcal E$ is simply connected. Indeed,
let $\beta$ be a simple loop in $\mathcal E^1$. Then it admits a
preimage $\ti\beta$ in $\Gamma$ which is either a loop; or a path
connecting two points $v_i{\in} A$ such that
$\s{pr}_E(v_i){=}v{\in} E\ (i{=}1,2).$ In the first case since
$C(\Gamma)$ is simply connected, the loop $\ti\beta$ is trivial
and so   $\beta$ is trivial in $\mathcal E$. In the second case we
have $v_i\in F_v$ for the set $F_v$ introduced in \ref{diric}. By
Proposition \ref{diricfund}, $F_v$ is $v$-star geodesic and there
exist two geodesics $l_i\subset F_v$ connecting $v_i$ with $v.$
The loop $\ti\eta{=}\ti\beta\cup l_1\cup l_2$ is trivial in
$C(\Gamma)$. We have $\s{pr}_E(F_v){=}\s{pr}_E(l_i){=}v$,  so
$\s{pr}_E(\ti\eta){=}\s{pr}_E(\ti\beta){=}\beta$ is as above a
trivial loop in $\mathcal E$. So $\mathcal E$ is simply connected.
Furthermore each relation in $H$ corresponds to a 2-disk $D$ in
$\mathcal E$ such that $D=\s{pr}_E(\ti D)$ where $\ti D$ is a
2-disk in $C(\Gamma).$ By construction the projection
$\s{pr}_E:C(\Gamma)^1\to \mathcal E^1$
 is an
isometric map. Therefore every relation in $H$ follows from
finitely many generating relations each of length at most $n$.
Thus the subgroup $H$ is finitely presented. \qed

\medskip

\bf Corollary.\sl\ Let a finitely presented group $G$ act
3-discontinuously and 2-cocompactly on a  compactum $T$. If   $H$
is a dynamically bounded subgroup of $G$ acting cocompactly on
$T\setminus \b\Lambda H$ then $H$ is finitely presented too.

\medskip

\bf Proof.\rm\ It follows immediately from Propositions
\ref{extens} and \ref{finpres}.

\medskip

\bf Remarks\rm. 1. In the above proof  we could at once assume
(w.l.o.g.) that $\Gamma{=}\Gamma_0$ is the Cayley graph. Indeed
the proper quasi-isometry $\varphi$ extends equivariantly to a
homeomorphism of $T$ keeping $\b\Lambda H$ invariant \cite[Lemma
2.5]{GP09}. This gives an equivariant proper map $\Gamma_0{\sqcup}
T\to \Gamma{\sqcup} T$ preserving $\b\Lambda H$. So the action of
$H$ on $(\Gamma_0{\sqcup} T)\setminus \b\Lambda H$ is cocompact
too.

2. The above Proposition was inspired by \cite[Theorem 1]{DG10}
establishing that the maximal parabolic subgroups of finitely
presented relatively hyperbolic groups are finitely presented.
This result follows from the above Corollary as maximal parabolic
subgroups act cocompactly outside their limit points and are
dynamically bounded (see Corollary of \ref{extens})

\medskip

\medskip

 \noindent We finish the Section by
a series of examples and questions.

\medskip

\vskip5pt \bf Examples\rm.\ 1)\ {\it An example of a subgroup
acting cocompactly on the complement of its limit set and which is
not a parabolic subgroup for any convergence action of  the
ambiant group}.

\medskip

 Let $G<{\rm Isom}\hn$ be a uniform
lattice. Let us fix two  elements $a$ and $b$ of $G$ having
different fixed points (i.e. generating a non-elementary subgroup
of $G$).

For a  sufficiently big $n_0$ the subgroup $H{=}<b,
a^{n_0}ba^{-n_0}>$  is free (Schottky) and quasiconvex in $G$.
Since $G$ contains no parabolics it follows that  the limit set of
$H$ is a proper Cantor subset of $\Bbb S^{n-1}$. The group $H$
acts geometrically finitely (without parabolic elements) on $\Bbb
S^{n-1}\setminus \b\Lambda H,$ and $({\Bbb S}^{n-1}\setminus \La
H)/H$ is a compact $(n{-}1)$-manifold homeomorphic to the
connected sum $(\Bbb S^{n-2}\times \Bbb S^1)\# (\Bbb S^{n-2}\times
\Bbb S^1)$. We have $gbg^{-1}{\in} H\cap gHg^{-1}$ where
$g{=}a^{n_0}{\not\in} H.$ Thus the subgroup $H{\cap}gHg^{-1}$ is
infinite. It is well known (e.g.  follows
  from our Proposition \ref{nearHoro}) that $H$ cannot be
parabolic for any geometrically finite action of $G$. \qed

\medskip

2) \ {\it An example of a dynamically bounded subgroup which is
not dynamically quasiconvex}.

\medskip

Take a 3-dimensional uniform arithmetic lattice $G<{\rm Isom}\Bbb
H^3$ such that $\Bbb H^3/G$ fibers over the circle. Let   $H$ be a
normal finitely generated subgroup of $G$ of infinite index which
is the group of the fiber manifold. It acts non-geometrically
finitely on $\Bbb H^n\ (n\geq 3).$ The group $G$ can be embedded
into another arithmetic lattice $G_0< {\rm Isom}\hn\ (n> 3).$
Since $G$ is (dynamically) quasiconvex in $G_0$, by Proposition
\ref{dcb} it is also dynamically bounded. Then by the hereditary
property (see Corollary of \ref{equivbound}) $H$ is a dynamically
bounded subgroup of $G_0.$

\medskip

3)\ {\it An example of a dynamically bounded subgroup of a
relatively hyperbolic finitely presented group which is not
finitely presented itself.}

\medskip

It is proved in \cite{KPV08} that any arithmetic non-uniform
lattice $G_0<{\rm Isom \Bbb H^n\ (n\geq 6)}$ contains a
geometrically finite subgroup $G<{\rm Isom} \Bbb H^4$ which
contains a normal finitely generated but infinitely presented
subgroup $F$. The subgroup $G$ is dynamically quasiconvex in
$G_0,$ and so is dynamically bounded. As in Example 2 by the
hereditary property $F$ is dynamically bounded too. This example
shows that Proposition \ref{finpres} is not true for dynamically
bounded subgroups without assuming the cocompactness of the action
outside   the limit set.

\medskip

Here are several questions which seem to be intriguing and open.

\vskip5pt \bf Questions\rm.\ 1) Suppose that the action $G\act\ti
T$ is 3-discontinuous and 2-cocompact. Let $H<G$ be a finitely
generated subgroup acting cocompactly on $T\setminus\b\Lambda H.$
Is it true that $H$ acts cocompactly on $\ti T\setminus\b\Lambda
H$ ?

\medskip

2) Suppose that the action $G\act\ti T$ is 3-discontinuous and
2-cocompact. Can $G$ contain a finitely generated subgroup $H$
such that $\La H\varsubsetneqq\La G$ and $H$ is not dynamically
bounded ?

In particular can a discrete finitely generated subgroup $H< {\rm
Isom}\Bbb H^n$ of a geometrically finite group $G$  such that $\La
H\varsubsetneqq\La G$ be not dynamically bounded ?

\medskip

\vskip5pt \bf Comments\rm.  If the  answer to the second question
is "no" then by \ref{extens} the answer to the first question is
  "yes". Then by \ref{subgrConvex} the subgroup $H$ acts cocompactly
  on $\ti T\setminus\b\Lambda H$ so
 $H$ is quasiconvex in $G.$ This would in particular imply that a
geometrically infinite finitely generated group in ${\rm Isom}
\hn$ acting cocompactly on the non-empty set $\Bbb
S^{n-1}\setminus \b\Lambda H$ ({\it totally degenerate} Kleinian
group) cannot appear as a subgroup of a lattice in ${\rm Isom}
\hn.$ This fact is true in dimension $n{=}3$ and follows from so
called {\it covering} theorem due to W. Thurston \cite[Proposition
7.1]{Mo84}. From the other hand totally degenerate Kleinian groups
exist in dimension $3$ and appear on the boundary of the
Teuchm\"uller spaces of surfaces \cite{Be70}. It is not known
whether they exist in higher dimensions.

We also note that every lattice $G$ in ${\rm Isom}\hn$ contains an
infinitely generated subgroup $H$ which is not dynamically
bounded. We thank Misha Kapovich for indicating to us that there
always exists  an infinitely generated subgroup of $G$ whose limit
set is a proper subset of $\Bbb S^{n-1}$ containing deep points
(see Introduction). This follows from the fact that the images of
$\b\Lambda H$ under a sequence of $(g_n){\subset} G\setminus H$
whose repeller point is a conical limit point of $G$, belonging to
$\b\Lambda H$, are dense in $\Bbb S^{n-1}$ [Ka00, Section 8.5].

\end{document}